# Laws of large numbers in stochastic geometry with statistical applications

MATHEW D. PENROSE

*Department of Mathematical Sciences, University of Bath, Bath BA2 7AY, United Kingdom.
E-mail: m.d.penrose@bath.ac.uk*

Given $n$ independent random marked $d$-vectors (points) $X_i$ distributed with a common density, define the measure $\nu_n = \sum_i \xi_i$, where $\xi_i$ is a measure (not necessarily a point measure) which stabilizes; this means that $\xi_i$ is determined by the (suitably rescaled) set of points near $X_i$. For bounded test functions $f$ on $R^d$, we give weak and strong laws of large numbers for $\nu_n(f)$. The general results are applied to demonstrate that an unknown set $A$ in $d$-space can be consistently estimated, given data on which of the points $X_i$ lie in $A$, by the corresponding union of Voronoi cells, answering a question raised by Khmaladze and Toronjadze. Further applications are given concerning the Gamma statistic for estimating the variance in nonparametric regression.

## 1. Introduction

Many interesting random variables in stochastic geometry arise as sums of contributions from each point of a point process $\mathcal{X}_n$ comprising $n$ independent random $d$-vectors $X_i, 1 \leq i \leq n$, distributed with common density function. General limit theorems, including laws of large numbers (LLNs), central limit theorems and large deviation principles, have been obtained for such variables, based on a notion of *stabilization* (local dependence) of the contributions; see [16, 17, 18, 20]. In particular, Penrose and Yukich [18] derive a general weak LLN of the form

$$\sum_{i=1}^n n^{-1} \xi(n^{1/d} X_i; \{n^{1/d} X_1, \ldots, n^{1/d} X_n\}) \xrightarrow{L^1} \beta, \qquad (1.1)$$

where $\xi(x; \mathcal{X})$ is a translation-invariant, real-valued functional defined for all finite $\mathcal{X} \subset \mathbb{R}^d$ and $x \in \mathcal{X}$, with $\xi$ satisfying stabilization and $(1+\varepsilon)$th moment conditions. The result of [18] also provides information about the limiting constant $\beta$.

Numerous applications of (1.1) are given in [18], for example, to sums of the form $\sum_{e \in G(\mathcal{X}_n)} \phi(n^{1/d}|e|)$, where $\phi$ is a test function satisfying a growth bound, $G(\mathcal{X}_n)$ is (for example) the nearest neighbour or Delaunay graph on $\mathcal{X}_n$ and $|e|$ denotes the length of edge $e$. For such sums, Jimenez and Yukich [9] give a *strong* LLN. Other examples considered in [2, 18] and elsewhere include those concerned with Voronoi graphs, minimal spanning trees, germ–grain models and on-line packing. In the present paper, we shall give general results extending the basic limit theorem (1.1) in the following directions.

*1. Almost sure convergence.* Under a moments condition and growth bound on the increment in the left-hand side of (1.1) on addition of a further point (the so-called 'add one cost'), this result can be extended to a strong LLN.

*2. Convergence of measures.* Consider the random measure on $\mathbb{R}^d$ comprising a point mass at each point $X_i$ of $\mathcal{X}_n$, of magnitude given by the $i$th term in the sum in the left-hand side of (1.1). This measure keeps track of the spatial locations of the contributions to the sum. Its asymptotic behaviour has been considered recently in [2, 7, 14, 20]. In fact, it is not necessary to restrict oneself to a point mass at $X_i$ and one can generalize further by considering the case where the contribution of the $i$th term to the random measure is some measure determined by $(X_i, \mathcal{X}_n)$ and localized 'near' to $X_i$ in some sense; for example, the 1-dimensional Hausdorff measure on the graph $G(\mathcal{X}_n)$ takes this form. We provide LLNs for the integrals of appropriate test functions against such measures.

*3. Non-translation-invariance.* It turns out that the translation-invariance condition on $\xi$ (which will be defined formally later) can also be relaxed in the limiting result (1.1) if, instead of scaling the point process $\mathcal{X}$ globally, as in (1.1), one scales the point process $\mathcal{X}_n$ locally at each point $X_i$ to keep the average density of points bounded, as $n$ becomes large.

*4. Marked point processes.* In a number of examples in the literature, the points of the point process $\mathcal{X}_n$ are required to carry independent identically distributed marks in some abstract measure space (the mark space).

We state our general results in Section 2 and prove them in Section 4, using auxiliary results on weak convergence given in Section 3. In Section 5 (which can be read without consulting Sections 3–4), we illustrate our general results with two specific applications in nonparametric statistics. Many other fields of application have been discussed elsewhere [2, 16, 18] and we do not attempt to review all of these, but our results often enable us to generalize results in those papers as indicated above.

The first application in Section 5 is to a question raised by Khmaladze and Toronjadze [11], motivated by statistical image analysis. Suppose $A$ is an unknown subset of the unit cube in $\mathbb{R}^d$ and that the random vectors $X_i$ represent 'sensors' which are scattered at random over the unit cube. Each sensor can detect whether it lies in the set $A$ or not. A reasonable estimate $A_n$ of the unknown set $A$, given the binary data from $n$ sensor locations $\mathcal{X}_n$, is then given by the union of the Voronoi cells, relative to the random point set $\mathcal{X}_n$, of those points $X_i$ which lie in $A$. Essentially, the question raised in [11] is whether or not $A_n$ is a consistent estimator for $A$ as $n \to \infty$; we answer this affirmatively via the strong LLN for measures.

Our second application is to a nonparametric regression model. Suppose that with each point $X_i \in \mathbb{R}^d$, we have an associated real-valued measurement $Y_i$, related to $X_i$ by $Y_i = h(X_i) + e_i$. Here, we assume that the function $h \in C^2(\mathbb{R}^d, \mathbb{R})$ is unknown and that the independent error $e_i$ has mean zero and unknown variance $\sigma^2$. For a large observed sample $(X_i, Y_i)_{i=1}^n$, a possible estimate for $\sigma^2$ is given by $W_i := (Y_{j(i,n,1)} - Y_i)/2$, where $j(i,n,1)$ is chosen so that $X_{j(i,n,1)}$ is the nearest neighbour of $X_i$ in $\mathcal{X}_n$. This estimator is approximately unbiased for $n$ large since then $W_i \approx (e_{j(i,n,1)} - e_i)^2/2$.

The mean $\overline{W}$ of $W_1, \ldots, W_n$ is the so-called *Gamma statistic* [5] based on nearest neighbours; one can consider a similar statistic based on $k$th nearest neighbours. We shall use our general results (requiring non-translation-invariance and marked points) to derive large-sample asymptotic properties associated with Gamma statistics.

## 2. Notation and general results

Let $(\mathcal{M}, \mathcal{F}_\mathcal{M}, \mu_\mathcal{M})$ be a probability space (the *mark space*). Let $\xi(\mathbf{x}; \mathcal{X}, A)$ be a Borel measurable $\mathbb{R}$-valued function defined for all triples $(\mathbf{x}, \mathcal{X}, A)$, where $\mathcal{X} \subset \mathbb{R}^d \times \mathcal{M}$ is finite and where $\mathbf{x} = (x, t) \in \mathcal{X}$ (so $x \in \mathbb{R}^d$ and $t \in \mathcal{M}$) and $A$ is a Borel set in $\mathbb{R}^d$. We assume that $\xi(\mathbf{x}; \mathcal{X}) := \xi(\mathbf{x}; \mathcal{X}, \cdot)$ is a $\sigma$-finite measure on $\mathbb{R}^d$.

Suppose $\mathbf{x} = (x, t) \in \mathbb{R}^d \times \mathcal{M}$ and $\mathcal{X} \subset \mathbb{R}^d \times \mathcal{M}$ is finite. If $\mathbf{x} \notin \mathcal{X}$, we abbreviate and write $\xi(\mathbf{x}; \mathcal{X})$ instead of $\xi(\mathbf{x}; \mathcal{X} \cup \{\mathbf{x}\})$. We also write $\mathcal{X}^\mathbf{x}$ for $\mathcal{X} \cup \{\mathbf{x}\}$. Given $y \in \mathbb{R}^d$ and $a \in \mathbb{R}$, we set $y + a\mathbf{x} := (y + ax, t)$. Let $y + a\mathcal{X} := \{(y + a\mathbf{w}) : \mathbf{w} \in \mathcal{X}\}$; in other words, scalar multiplication and translation act only on the first component of elements of $\mathbb{R}^d \times \mathcal{M}$. For $A \subseteq \mathbb{R}^d$, we set $y + aA = \{Y + aw : w \in A\}$. We say $\xi$ is *translation invariant* if

$$\xi(\mathbf{x}; \mathcal{X}, A) = \xi(y + \mathbf{x}; y + \mathcal{X}, y + A)$$

for all $y \in \mathbb{R}^d$, all finite $\mathcal{X} \subset \mathbb{R}^d \times \mathcal{M}$ and $\mathbf{x} \in \mathcal{X}$ and all Borel $A \subseteq \mathbb{R}^d$.

Let $\kappa$ be a probability density function on $\mathbb{R}^d$. Abusing notation slightly, we also let $\kappa$ denote the corresponding probability measure on $\mathbb{R}^d$, that is, we write $\kappa(A)$ for $\int_A \kappa(x) \, dx$, for Borel $A \subseteq \mathbb{R}^d$. For all $\lambda > 0$, let $\lambda\kappa$ denote the measure on $\mathbb{R}^d$ with density $\lambda\kappa(\cdot)$ and let $\mathcal{P}_\lambda$ denote a Poisson point process in $\mathbb{R}^d \times \mathcal{M}$ with intensity measure $\lambda\kappa \times \mu_\mathcal{M}$.

Let $(X, T), (X', T'), (X_1, T_1), (X_2, T_2), \ldots$ denote a sequence of independent identically distributed random elements of $\mathbb{R}^d \times \mathcal{M}$ with distribution $\kappa \times \mu_\mathcal{M}$ and set $\mathbf{X} := (X, T)$, $\mathbf{X}' := (X', T')$ and $\mathbf{X}_i := (X_i, T_i)$, $i \geq 1$. For $n \in \mathbb{N}$, let $\mathcal{X}_n$ be the point process in $\mathbb{R}^d \times \mathcal{M}$ given by $\mathcal{X}_n := \{\mathbf{X}_1, \mathbf{X}_2, \ldots, \mathbf{X}_n\}$. Let $\mathcal{H}_\lambda$ denote a Poisson point process in $\mathbb{R}^d \times \mathcal{M}$ with intensity $\lambda$ times the product of $d$-dimensional Lebesgue measure and $\mu_\mathcal{M}$ (i.e., a homogeneous marked Poisson process in $\mathbb{R}^d$ with intensity $\lambda$) and let $\tilde{\mathcal{H}}_\lambda$ denote an independent copy of $\mathcal{H}_\lambda$.

Suppose we are given a family of non-empty open subsets $\Omega_\lambda$ of $\mathbb{R}^d$, indexed by $\lambda \geq 1$, that are non-decreasing in $\lambda$, that is, satisfying $\Omega_\lambda \subseteq \Omega_{\lambda'}$ for $\lambda < \lambda'$. Denote by $\Omega_\infty$ the limiting set, that is, set $\Omega_\infty := \bigcup_{\lambda \geq 1} \Omega_\lambda$. Suppose we are given a further Borel set $\Omega$ with $\Omega_\infty \subseteq \Omega \subseteq \mathbb{R}^d$. In many examples, one takes $\Omega_\lambda = \Omega$ for all $\lambda$, either with $\Omega = \mathbb{R}^d$ or with $\kappa$ supported by $\Omega$.

For $\lambda > 0$ and for finite $\mathcal{X} \subset \mathbb{R}^d \times \mathcal{M}$ with $\mathbf{x} = (x, t) \in \mathcal{X}$ and Borel $A \subset \mathbb{R}^d$, let

$$\xi_\lambda(\mathbf{x}; \mathcal{X}, A) := \xi(\mathbf{x}; x + \lambda^{1/d}(-x + \mathcal{X}), x + \lambda^{1/d}(-x + A))\mathbf{1}_{\Omega_\lambda}(x). \tag{2.1}$$

When $\xi$ is translation invariant, the rescaled measure $\xi_\lambda$ simplifies to

$$\xi_\lambda(\mathbf{x}; \mathcal{X}, A) = \xi(\lambda^{1/d}\mathbf{x}; \lambda^{1/d}\mathcal{X}, \lambda^{1/d}A)\mathbf{1}_{\Omega_\lambda \times \mathcal{M}}(\mathbf{x}). \tag{2.2}$$

In general, the point process $x + \lambda^{1/d}(-x + \mathcal{X})$ is obtained by a dilation, centred at $x$, of the original point process. Loosely speaking, this dilation has the effect of reducing the density of points by a factor of $\lambda$. Thus, for $\mathbf{x} = (x,t) \in \mathbb{R}^d \times \mathcal{M}$, the rescaled measure $\xi_\lambda(\mathbf{x}; \mathcal{X}, A)$ is the original measure $\xi$ at $x$ relative to the image of the point process $\mathcal{X}$ under a dilation about $x$, acting on the image of 'space' (i.e., the set $A$) under the same dilation. This 'dilation of space' has the effect of concentrating the measure near to $x$; for example, if $\xi(\mathbf{x}; \mathcal{X})$ is a unit point mass at $x + y(x)$ for some measurable choice of function $x \mapsto y(x) \in \mathbb{R}^d$ and $\Omega_\lambda = \mathbb{R}^d$, then $\xi_\lambda(x,t; \mathcal{X})$ would be a unit point mass at $x + \lambda^{-1/d}y(x)$.

Our principal objects of interest are the random measures $\nu_{\lambda,n}$ on $\mathbb{R}^d$, defined for $\lambda > 0$ and $n \in \mathbb{N}$ by $\nu_{\lambda,n} := \sum_{i=1}^n \xi_\lambda(\mathbf{X}_i; \mathcal{X}_n)$. We study these measures via their action on test functions in the space $B(\Omega)$ of bounded Borel measurable functions on $\Omega$. We let $\tilde{B}(\Omega)$ denote the subclass of $B(\Omega)$ consisting of those functions that are Lebesgue-almost everywhere continuous. When $\Omega \neq \mathbb{R}^d$, we extend functions $f \in B(\Omega)$ to $\mathbb{R}^d$ by setting $f(x) = 0$ for $x \in \mathbb{R}^d \setminus \Omega$. Given $f \in B(\Omega)$, set $\langle f, \xi_\lambda(\mathbf{x}; \mathcal{X}) \rangle := \int_{\mathbb{R}^d} f(z)\xi_\lambda(\mathbf{x}; \mathcal{X}, dz)$. Also, set

$$\langle f, \nu_{\lambda,n} \rangle := \int_\Omega f \, d\nu_{\lambda,n} = \sum_{i=1}^n \langle f, \xi_\lambda(\mathbf{X}_i; \mathcal{X}_n) \rangle. \tag{2.3}$$

The indicator function $\mathbf{1}_{\Omega_\lambda}(x)$ in the definition (2.1) of $\xi_\lambda$ means that only points $\mathbf{X}_i \in \Omega_\lambda \times \mathcal{M}$ contribute to $\nu_{\lambda,n}$. In most examples, the sets $\Omega_\lambda$ are all the same and often are all $\mathbb{R}^d$. However, there are cases where moment conditions such as (2.5) below hold for a sequence of sets $\Omega_\lambda$, but would not hold if we were to take $\Omega_\lambda = \Omega$ for all $\lambda$; see, for example, [15]. Likewise, in some examples, the measure $\xi(\mathbf{x}; \mathcal{X})$ is not finite on the whole of $\mathbb{R}^d$, but is well behaved on $\Omega$, hence the restriction of attention to test functions in $B(\Omega)$.

Let $|\cdot|$ denote the Euclidean norm on $\mathbb{R}^d$ and for $x \in \mathbb{R}^d$ and $r > 0$, define the ball $B_r(x) := \{y \in \mathbb{R}^d : |y - x| \leq r\}$. We denote by $\mathbf{0}$ the origin of $\mathbb{R}^d$ and abbreviate $B_r(\mathbf{0})$ to $B_r$. We write $B_r^*(x)$ for $B_r(x) \times \mathcal{M}$, $B_r^*$ for $B_r \times \mathcal{M}$ and $(B_r^*)^c$ for $(\mathbb{R}^d \setminus B_r) \times \mathcal{M}$. Finally, we let $\omega_d$ denote the Lebesgue measure of the $d$-dimensional unit ball $B_1$.

We say a set $\mathcal{X} \subset \mathbb{R}^d \times \mathcal{M}$ is *locally finite* if $\mathcal{X} \cap B_r^*$ is finite for all finite $r$. For $\mathbf{x} \in \mathbb{R}^d \times \mathcal{M}$ and Borel $A \subseteq \mathbb{R}^d$, we extend the definition of $\xi(\mathbf{x}; \mathcal{X}, A)$ to locally finite infinite point sets $\mathcal{X}$ by setting

$$\xi(\mathbf{x}; \mathcal{X}, A) := \limsup_{K \to \infty} \xi(\mathbf{x}; \mathcal{X} \cap B_K^*, A).$$

Also, for $\mathbf{x} = (x,t) \in \mathbb{R}^d \times \mathcal{M}$, we define the $x$-shifted version $\xi_\infty^\mathbf{x}(\cdot,\cdot)$ of $\xi(\mathbf{x}; \cdot, \cdot)$ by

$$\xi_\infty^\mathbf{x}(\mathcal{X}, A) = \xi(\mathbf{x}; x + \mathcal{X}, x + A).$$

Note that if $\xi$ is translation invariant, then $\xi_\infty^{(x,t)}(\mathcal{X}, A) = \xi_\infty^{(0,t)}(\mathcal{X}, A)$ for all $x \in \mathbb{R}^d$, $t \in \mathcal{M}$ and Borel $A \subseteq \mathbb{R}^d$.

The following notion of *stabilization* is similar to those used in [2, 18].

**Definition 2.1.** *For any locally finite $\mathcal{X} \subset \mathbb{R}^d \times \mathcal{M}$ and any $\mathbf{x} = (x,t) \in \mathbb{R}^d \times \mathcal{M}$, define $R(\mathbf{x}; \mathcal{X})$ (the radius of stabilization of $\xi$ at $\mathbf{x}$ with respect to $\mathcal{X}$) to be the smallest integer-valued $r$ such that $r \geq 0$ and*

$$\xi(\mathbf{x}; x + [\mathcal{X} \cap B_r^*] \cup \mathcal{Y}, B) = \xi(\mathbf{x}; x + [\mathcal{X} \cap B_r^*], B)$$

*for all finite $\mathcal{Y} \subseteq (B_r^*)^c$ and Borel $B \subseteq \mathbb{R}^d$. If no such $r$ exists, we set $R(\mathbf{x}; \mathcal{X}) = \infty$.*

In the case where $\xi$ is translation invariant, $R((x,t); \mathcal{X}) = R((\mathbf{0}, t); \mathcal{X})$, so $R((x,t); \mathcal{X})$ does not depend on $x$. Of particular importance to us will be radii of stabilization with respect to the homogeneous Poisson processes $\mathcal{H}_\lambda$.

We assert that $R(\mathbf{x}; \mathcal{X})$ is a measurable function of $\mathcal{X}$ and hence, when $\mathcal{X}$ is a random point set such as $\mathcal{H}_\lambda$, $R(\mathbf{x}; \mathcal{X})$ is an $\mathbb{N} \cup \{\infty\}$-valued random variable. To see this assertion, observe that by Dynkin's pi-lambda lemma, for any $k \in \mathbb{N}$, the event $\{R(\mathbf{x}; \mathcal{X}) \leq k\}$ equals the event $\bigcap_{B \in \mathcal{B}} \{s(\mathcal{X}, B) = i(\mathcal{X}, B)\}$, where $\mathcal{B}$ is the $\Pi$-system consisting of the rectilinear hypercubes in $\mathbb{R}^d$ whose corners have rational coordinates and for $B \in \mathcal{B}$, we set

$$s(\mathcal{X}, B) := \sup\{\xi(\mathbf{x}; x + ([\mathcal{X} \cap B_k^*] \cup \mathcal{Y}), B) : \mathcal{Y} \subseteq (B_k^*)^c\},$$
$$i(\mathcal{X}, B) := \inf\{\xi(\mathbf{x}; x + ([\mathcal{X} \cap B_k^*] \cup \mathcal{Y}), B) : \mathcal{Y} \subseteq (B_k^*)^c\}.$$

Also, $s(\mathcal{X}, B)$ is a measurable function of $\mathcal{X}$ because we assume $\xi$ is Borel measurable and, for any $b$, we have

$$\{\mathcal{X} : s(\mathcal{X}) > b\} = \pi_1(\{(\mathcal{X}, \mathcal{Y}) : \xi(\mathbf{x}; x + [\mathcal{X} \cap B_k^*] \cup [\mathcal{Y} \setminus B_k^*], B) > b\}),$$

where $\pi_1$ denotes projection onto the first component, acting on pairs $(\mathcal{X}, \mathcal{Y})$, with $\mathcal{X}$ and $\mathcal{Y}$ finite sets in $\mathbb{R}^d \times \mathcal{M}$. Similarly, $i(\mathcal{X}, B)$ is a measurable function of $\mathcal{X}$.

For $\mathbf{x} = (x,t) \in \mathbb{R}^d \times \mathcal{M}$, let $\xi_\lambda^*(\mathbf{x}; \mathcal{X}, \cdot)$ be the point measure at $x$ with total measure $\xi_\lambda(\mathbf{x}; \mathcal{X}, \Omega)$, that is, for Borel $A \subseteq \mathbb{R}^d$, let

$$\xi_\lambda^*(\mathbf{x}; \mathcal{X}, A) := \xi_\lambda(\mathbf{x}; \mathcal{X}, \Omega) \mathbf{1}_A(x). \tag{2.4}$$

We consider measures $\xi$ and test functions $f \in B(\Omega)$ satisfying one of the following assumptions.

A1: $\xi((x,t); \mathcal{X}, \cdot)$ is a point mass at $x$ for all $(x, t, \mathcal{X})$.
A2: $\xi(\mathbf{x}; \mathcal{X}, \cdot)$ is absolutely continuous with respect to Lebesgue measure on $\mathbb{R}^d$, with Radon–Nikodym derivative denoted $\xi'(\mathbf{x}; \mathcal{X}, y)$ for $y \in \mathbb{R}^d$, satisfying $\xi'(\mathbf{x}; \mathcal{X}, y) \leq K_0$ for all $(\mathbf{x}; \mathcal{X}, y)$, where $K_0$ is a finite positive constant.
A3: $f$ is almost everywhere continuous, that is, $f \in \tilde{B}(\Omega)$.

Note that assumption A1 implies that $\xi_\lambda^* = \xi_\lambda$, and that assumption A2 will hold if $\xi(\mathbf{x}, \mathcal{X}, \cdot)$ is Lebesgue measure on some random subset of $\mathbb{R}^d$ determined by $\mathbf{x}, \mathcal{X}$.

Our first general result is a weak law of large numbers for $\langle f, \nu_{\lambda,n} \rangle$ defined at (2.3), for $f \in B(\Omega)$. This extends [18], which is concerned only with the case where $f$ is a constant. We require almost surely finite radii of stabilization with respect to homogeneous Poisson processes, along with a moments condition.

**Theorem 2.1.** *Suppose that $R((x,T); \mathcal{H}_{\kappa(x)})$ is almost surely finite for $\kappa$-almost all $x \in \Omega_\infty$. Suppose, also, that $f \in B(\Omega)$ and that one or more of assumptions* A1, A2, A3 *holds. Let $q = 1$ or $q = 2$ and let the sequence $(\lambda(n), n \geq 1)$ be a sequence of positive numbers with $\lambda(n)/n \to 1$ as $n \to \infty$. If there exists $p > q$ such that*

$$\limsup_{n \to \infty} \mathbb{E}[\xi_{\lambda(n)}(\mathbf{X}; \mathcal{X}_{n-1}, \Omega)^p] < \infty, \tag{2.5}$$

*then, as $n \to \infty$, we have the $L^q$ convergence*

$$n^{-1} \langle f, \nu_{\lambda(n), n} \rangle \to \int_{\Omega_\infty} f(x) \mathbb{E}[\xi_\infty^{(x,T)}(\mathcal{H}_{\kappa(x)}, \mathbb{R}^d)] \kappa(x) \, dx \qquad as \ n \to \infty, \tag{2.6}$$

*with finite limit, and the $L^1$ convergence*

$$n^{-1} \sum_{i=1}^n |\langle f, \xi_{\lambda(n)}(\mathbf{X}_i; \mathcal{X}_n) - \xi_{\lambda(n)}^*(\mathbf{X}_i; \mathcal{X}_n) \rangle| \to 0 \qquad as \ n \to \infty. \tag{2.7}$$

To extend Theorem 2.1 to a strong law, we need to assume extra conditions concerning the so-called *add one cost*, that is, the effect of adding a single further point on the measure $\nu_{\lambda(n), n-1}$. We define three different types of add one cost, the first two of which refer to a test function $f$. Given $f \in B(\Omega)$ and $n \geq 2$, $\lambda \geq 1$, set $\Delta_{\lambda,n}(f) := \langle f, \nu_{\lambda,n} - \nu_{\lambda,n-1} \rangle$. Also, let $G_{\lambda,n}(f)$ be the sum in (2.7), that is, set

$$G_{\lambda,n}(f) := \sum_{i=1}^n |\langle f, \xi_\lambda(\mathbf{X}_i; \mathcal{X}_n) - \xi_\lambda^*(\mathbf{X}_i; \mathcal{X}_n) \rangle|$$

and set $\Delta'_{\lambda,n}(f) := G_{\lambda,n}(f) - G_{\lambda,n-1}(f)$. Furthermore, for $1 \leq i \leq n-1$, let $|\xi_\lambda(\mathbf{X}_i; \mathcal{X}_n) - \xi_\lambda(\mathbf{X}_i; \mathcal{X}_{n-1})|(\Omega)$ denote the total variation (i.e., the sum of positive and negative parts) of the signed measure $\xi_\lambda(\mathbf{X}_i; \mathcal{X}_n, \cdot) - \xi_\lambda(\mathbf{X}_i; \mathcal{X}_{n-1}, \cdot)$ on $\Omega$ and define $\tilde{\Delta}_{\lambda,n}$ by

$$\tilde{\Delta}_{\lambda,n} := \xi_\lambda(\mathbf{X}_n; \mathcal{X}_n, \Omega) + \sum_{i=1}^{n-1} |\xi_\lambda(\mathbf{X}_i; \mathcal{X}_n) - \xi_\lambda(\mathbf{X}_i; \mathcal{X}_{n-1})|(\Omega). \tag{2.8}$$

Given a random variable $W$, as usual, we let $\|W\|_p := \mathbb{E}[|W|^p]^{1/p}$ for $p = 1$ and $\|W\|_\infty := \inf\{t > 0 : P[|W| > t] = 0\}$ with $\inf(\varnothing) := +\infty$.

**Theorem 2.2.** *Suppose that $R((x,T); \mathcal{H}_{\kappa(x)})$ is almost surely finite for $\kappa$-almost all $x \in \Omega_\infty$. Suppose, also, that $f \in B(\Omega)$ and that one or more of assumptions* A1, A2, A3 *holds. Suppose $\lambda(n)/n \to 1$ as $n \to \infty$ and that there exists $p > 1$ such that (2.5) holds. Suppose $\beta \geq 1$ and $p' > 2(\beta + 1)$. If, as $n \to \infty$, we have*

$$\|\Delta_{\lambda(n),n}(f)\|_\infty = O(n^\beta), \qquad \|\Delta_{\lambda(n),n}(f)\|_{p'} = O(1), \tag{2.9}$$

*then (2.6) holds with almost sure convergence. If, instead of (2.9), we have*

$$\|\Delta'_{\lambda(n),n}(f)\|_\infty = O(n^\beta), \qquad \|\Delta'_{\lambda(n),n}(f)\|_{p'} = O(1), \tag{2.10}$$

*then (2.7) holds with almost sure convergence. Finally, if*

$$\|\tilde{\Delta}_{\lambda(n),n}\|_\infty = O(n^\beta), \qquad \|\tilde{\Delta}_{\lambda(n),n}\|_{p'} = O(1), \tag{2.11}$$

*then both (2.9) and (2.10) hold, so both (2.6) and (2.7) hold with almost sure convergence.*

***Remarks.*** Certain weak laws of large numbers for $\langle f, \nu^\xi_{\lambda(n),n}\rangle$ follow directly from Theorem 2.4 of Baryshnikov and Yukich [2]. However, the conditions in Theorem 2.1 are weaker in many ways than those in [2], as one might expect, since we consider only the law of large numbers, whereas [2] is concerned with Gaussian limits.

For example, in [2], attention is restricted to cases where assumption A1 holds. It is often natural to drop this restriction. Also, in [2], it is assumed that $\kappa$ has compact convex support and is continuous on its support, whereas we make no assumptions here on $\kappa$. Moreover, in [2], attention is restricted to *continuous* bounded test functions $f$, whereas we consider test functions which are merely bounded (under A1 or A2) or bounded and almost everywhere continuous (under A3). Thus, we can consider test functions which are indicator functions of Borel sets $A$ in $\Omega$.

Our stabilization conditions refer only to homogeneous Poisson processes, and not to any non-homogeneous Poisson processes, as in [2]; unlike that paper, we require only that radii of stabilization be almost surely finite, with no condition on their tails. Also, our moments condition (2.5) is simpler than the corresponding condition in [2] (equation (2.2) of [2]).

Almost sure convergence, that is, the strong law of large numbers, is not addressed in [2] or [18]. Some strong laws for graphs arising in geometric probability are derived by Jimenez and Yukich [9] and we add to these. As in [9], we actually prove *complete convergence*, as defined in, for example, [13] or [17].

Unlike [2, 9, 18], we spell out the statement and proof of our law of large numbers for *marked* point processes (i.e., point processes in $\mathbb{R}^d \times \mathcal{M}$, rather than in $\mathbb{R}^d$). This setting includes many interesting examples, such as germ–grain models and on-line packing, and generalizes the unmarked point process setting because we can always take $\mathcal{M}$ to have a single element and then identify $\mathbb{R}^d \times \mathcal{M}$ with $\mathbb{R}^d$ to recover results for unmarked point processes from the general results for marked point processes.

*Poisson samples.* It is also of interest to obtain a similar result to Theorem 2.1 for the random measure $\mu_\lambda$ defined in a similar manner to $\nu_{\lambda,n}$, but using the Poisson point process $\mathcal{P}_\lambda$ instead of $\mathcal{X}_n$, taking $\lambda \to \infty$. Such a result can indeed be obtained by a similar proof, with $L^1$ convergence if, instead of the moments condition (2.5), one assumes

$$\sup_{\lambda \geq 1} \mathbb{E}[\xi_\lambda(\mathbf{X}; \mathcal{P}_\lambda, \Omega)^p] < \infty, \qquad \text{for some } p > 1.$$

*Multisample statistics.* Suppose that $\mathcal{X}_{n_1} \subset \mathbb{R}^d$ represents a sample of $n_1$ points of 'type one' and $\mathcal{Y}_{n_2} \subset \mathbb{R}^d$ represents an independent sample of $n_2$ points of 'type two,' possibly having a different underlying density function. Suppose that for $i = 1, 2$, the functionals $\xi^{(1)}(x; \mathcal{X}, \mathcal{Y})$ and $\xi^{(2)}(y; \mathcal{X}, \mathcal{Y})$ are defined in a translation-invariant and stabilizing manner for finite $\mathcal{X} \subset \mathbb{R}^d$, finite $\mathcal{Y} \subset \mathbb{R}^d$, $x \in \mathcal{X}$ and $y \in \mathcal{Y}$. When the integers $n_i = n_i(n)$ ($i = 1, 2$) satisfy $n_1(n) + n_2(n) = n$ and $n_i(n)/n \to \pi_i \in (0,1)$ as $n \to \infty$, the sum

$$\sum_{i=1}^{n_1} \xi^{(1)}(n^{1/d} X_i; n_1^{1/d} \mathcal{X}_{n_1}, n^{1/d} \mathcal{Y}_{n_2}) + \sum_{j=1}^{n_2} \xi^{(2)}(n^{1/d} Y_j; n^{1/d} \mathcal{X}_{n_1}, n^{1/d} \mathcal{Y}_{n_2})$$

satisfies an LLN under weaker (or at least, different) conditions than those for the main result of Henze and Voigt ([8], Theorem 2.3) for such sums (and likewise for more than two samples). The proof is based on the techniques of this paper and the results are comparable to those for independently marked points where $\mathcal{M} = \{1, 2\}$ and the 'mark' determines whether a point is in sample $\mathcal{X}$ or $\mathcal{Y}$. We omit the details here.

## 3. Weak convergence and the objective method

In this section, we derive certain weak convergence results (Lemmas 3.3–3.6). We use a version of the 'objective method' [1, 21], whereby convergence in distribution (denoted $\xrightarrow{\mathcal{D}}$) for a functional defined on a sequence of finite probabilistic objects (in this case, rescaled marked point processes) is established by showing that these probabilistic objects themselves converge in distribution to an infinite probabilistic object (in this case, a homogeneous marked Poisson process) and that the functional of interest is continuous.

A *point process* in $\mathbb{R}^d \times \mathcal{M}$ is an $\mathcal{L}$-valued random variable, where $\mathcal{L}$ denotes the space of locally finite subsets of $\mathbb{R}^d \times \mathcal{M}$. We use the following metric on $\mathcal{L}$:

$$D(\mathcal{A}, \mathcal{A}') = (\max\{K \in \mathbb{N} : \mathcal{A} \cap B_K^* = \mathcal{A}' \cap B_K^*\})^{-1}. \tag{3.1}$$

With this metric, $\mathcal{L}$ is a metric space which is complete but not separable. In the unmarked case where $\mathcal{M}$ has a single element, our choice of metric is *not* the same as the metric used in Section 5.3 of [21]. Indeed, for one-point unmarked sets, our metric generates the discrete topology rather than the Euclidean topology.

Recall (see, e.g., [13], [19]) that $x \in \mathbb{R}^d$ is a *Lebesgue point* of $\kappa$ if $\varepsilon^{-d} \int_{B_\varepsilon(x)} |\kappa(y) - \kappa(x)| \, dy$ tends to zero as $\varepsilon \downarrow 0$ and that the Lebesgue Density Theorem tells us that

almost every $x \in \mathbb{R}^d$ is a Lebesgue point of $\kappa$. For subsequent results, it is useful to define the region

$$\Omega_0 := \{x \in \Omega_\infty : \kappa(x) > 0, x \text{ a Lebesgue point of } \kappa(\cdot)\}. \tag{3.2}$$

**Lemma 3.1.** *Suppose $x \in \Omega_0$ and suppose $(y(\lambda), \lambda > 0)$ is an $\mathbb{R}^d$-valued function with $|y(\lambda) - x| = O(\lambda^{-1/d})$ as $\lambda \to \infty$. Then there exist coupled realizations $\mathcal{P}'_\lambda$ and $\mathcal{H}'_{\kappa(x)}$ of $\mathcal{P}_\lambda$ and $\mathcal{H}_{\kappa(x)}$, respectively, such that*

$$D(\lambda^{1/d}(-y(\lambda) + \mathcal{P}'_\lambda), \mathcal{H}'_{\kappa(x)}) \xrightarrow{P} 0 \qquad \text{as } \lambda \to \infty. \tag{3.3}$$

**Proof.** Let $\mathcal{H}^+$ denote a homogeneous Poisson process of unit intensity in $\mathbb{R}^d \times \mathcal{M} \times [0, \infty)$. Let $\mathcal{P}'_\lambda$ denote the image of the restriction of $\mathcal{H}^+$ to the set $\{(w, t, s) \in \mathbb{R}^d \times \mathcal{M} \times [0, \infty) : s \leq \lambda \kappa(w)\}$ under the mapping $(w, t, s) \mapsto (w, t)$. Let $\mathcal{H}'_{\kappa(x)}$, denote the image of the restriction of $\mathcal{H}^+$ to the set $\{(w, t, s) \in \mathbb{R}^d \times \mathcal{M} \times [0, \infty) : s \leq \lambda \kappa(x)\}$ under the mapping

$$(w, t, s) \mapsto (\lambda^{1/d}(w - y(\lambda)), s).$$

By the Mapping Theorem [12], $\mathcal{P}'_\lambda$ has the same distribution as $\mathcal{P}_\lambda$, while $\mathcal{H}'_{\kappa(x)}$ has the same distribution as $\mathcal{H}_{\kappa(x)}$.

The number of points of the point set

$$(\lambda^{1/d}(-y(\lambda) + \mathcal{P}'_\lambda) \triangle \mathcal{H}'_{\kappa(x)}) \cap B^*_K$$

equals the number of points $(X, T, S)$ of $\mathcal{H}^+$ with $X \in B_{\lambda^{-1/d}K}(y(\lambda))$ and with either $\lambda \kappa(x) < S \leq \lambda \kappa(X)$ or $\lambda \kappa(X) < S \leq \lambda \kappa(x)$. This is Poisson distributed with mean

$$\lambda \int_{B_{\lambda^{-1/d}K}(y(\lambda))} |\kappa(z) - \kappa(x)| \, dz,$$

which tends to zero because $x$ is assumed to be a Lebesgue point of $\kappa$, and (3.3) follows. $\square$

In the rest of this section, given $x \in \mathbb{R}^d$, we write $\mathbf{x}$ for $(x, T)$ (i.e., for the point $x$ equipped with a generic random mark $T$) and, given $y \in \mathbb{R}^d$, we write $\mathbf{y}$ for $(y, T')$.

**Lemma 3.2.** *Suppose $(x, y) \in \Omega_0 \times \Omega_0$ with $x \neq y$. Let $(\lambda(k), \ell(k), m(k))_{k \in \mathbb{N}}$ be a $((0, \infty) \times \mathbb{N} \times \mathbb{N})$-valued sequence satisfying $\lambda(k) \to \infty$, and $\ell(k)/\lambda(k) \to 1$ and $m(k)/\lambda(k) \to 1$ as $k \to \infty$. Then, as $k \to \infty$,*

$$(\lambda(k)^{1/d}(-x + \mathcal{X}_{\ell(k)}), \lambda(k)^{1/d}(-x + \mathcal{X}_{m(k)}), \lambda(k)^{1/d}(-y + \mathcal{X}_{m(k)}),$$
$$\lambda(k)^{1/d}(-x + \mathcal{X}^{\mathbf{y}}_{\ell(k)}), \lambda(k)^{1/d}(-x + \mathcal{X}^{\mathbf{y}}_{m(k)}), \lambda(k)^{1/d}(-y + \mathcal{X}^{\mathbf{x}}_{m(k)}))$$
$$\xrightarrow{\mathcal{D}} (\mathcal{H}_{\kappa(x)}, \mathcal{H}_{\kappa(x)}, \tilde{\mathcal{H}}_{\kappa(y)}, \mathcal{H}_{\kappa(x)}, \mathcal{H}_{\kappa(x)}, \tilde{\mathcal{H}}_{\kappa(y)}). \tag{3.4}$$

**Proof.** In this proof, we write simply $\lambda$ for $\lambda(k)$, $\ell$ for $\ell(k)$ and $m$ for $m(k)$. We use the following coupling. Suppose we are given $\lambda$. On a suitable probability space, let $\mathcal{P}$ and $\tilde{\mathcal{P}}$ be independent copies of $\mathcal{P}_\lambda$, independent of $\mathbf{X}_1, \mathbf{X}_2, \ldots$.

Let $\mathcal{P}'$ be the point process in $\mathbb{R}^d \times \mathcal{M}$ consisting of those points $(V, T) \in \mathcal{P}$ such that $|V - x| < |V - y|$, together with those points $(V', U') \in \tilde{\mathcal{P}}$ with $|V' - y| < |V' - x|$. Clearly $\mathcal{P}'$ is a Poisson process of intensity $\lambda \kappa \times \mu_\mathcal{M}$ on $\mathbb{R}^d \times \mathcal{M}$.

Let $\mathcal{H}'_{\kappa(x)}$ and $\tilde{\mathcal{H}}'_{\kappa(y)}$ be independent copies of $\mathcal{H}_{\kappa(x)}$ and $\mathcal{H}_{\kappa(y)}$, respectively. Assume $\mathcal{H}'_{\kappa(x)}$ and $\tilde{\mathcal{H}}'_{\kappa(y)}$ are independent of $(\mathbf{X}_1, \mathbf{X}_2, \mathbf{X}_3, \ldots)$. Using Lemma 3.1, assume, also, that $\mathcal{H}'_{\kappa(x)}$ is coupled to $\mathcal{P}$ and $\tilde{\mathcal{H}}'_{\kappa(y)}$ is coupled to $\tilde{\mathcal{P}}$ in such a way that as $k \to \infty$,

$$\max(D(\mathcal{H}'_{\kappa(x)}, \lambda^{1/d}(-x + \mathcal{P})), D(\tilde{\mathcal{H}}'_{\kappa(y)}, \lambda^{1/d}(-y + \tilde{\mathcal{P}}))) \xrightarrow{P} 0. \tag{3.5}$$

Let $N$ denote the number of points of $\mathcal{P}'$ (a Poisson variable with mean $\lambda$). Choose an ordering on the points of $\mathcal{P}'$, uniformly at random from all $N!$ possible such orderings. Use this ordering to list the points of $\mathcal{P}'$ as $\mathbf{W}_1, \mathbf{W}_2, \ldots, \mathbf{W}_N$. Also, set $\mathbf{W}_{N+1} = \mathbf{X}_1, \mathbf{W}_{N+2} = \mathbf{X}_2, \mathbf{W}_{N+3} = \mathbf{X}_3$ and so on. Set

$$\mathcal{X}'_\ell := \{\mathbf{W}_1, \ldots, \mathbf{W}_\ell\}, \qquad \mathcal{X}'_m := \{\mathbf{W}_1, \ldots, \mathbf{W}_m\}.$$

Then $(\mathcal{X}'_\ell, \mathcal{X}'_m) \stackrel{\mathcal{D}}{=} (\mathcal{X}_\ell, \mathcal{X}_m)$ and $(\mathcal{H}'_{\kappa(x)}, \tilde{\mathcal{H}}'_{\kappa(y)}) \stackrel{\mathcal{D}}{=} (\mathcal{H}_{\kappa(x)}, \tilde{\mathcal{H}}_{\kappa(y)})$, where $\stackrel{\mathcal{D}}{=}$ denotes equality of distribution.

Let $K \in \mathbb{N}$ and let $\delta > 0$. Define the events

$$E := \{\mathcal{X}'_m \cap B^*_{\lambda^{-1/d}K}(x) = \mathcal{P}' \cap B^*_{\lambda^{-1/d}K}(x)\},$$
$$F := \{(\lambda^{1/d}(-x + \mathcal{P})) \cap B^*_K = \mathcal{H}'_{\kappa(x)} \cap B^*_K\}.$$

Event $E$ occurs unless either one or more of the $(N - m)^+$ "discarded" points of $\mathcal{P}'$ or one or more of the $(m - N)^+$ "added" points of $\{\mathbf{X}_1, \mathbf{X}_2, \ldots\}$ lies in $B^*_{\lambda^{-1/d}K}(x)$. For each added or discarded point, for sufficiently large $\lambda$, the probability of lying in $B^*_{\lambda^{-1/d}K}(x)$ is at most $\omega_d(\kappa(x) + 1)K^d/\lambda$ because $x$ is a Lebesgue point. Thus, for $k$ sufficiently large that $|m - \lambda| \leq \delta\lambda$, we have

$$P[E^c] \leq P[|N - \lambda| > \delta\lambda] + (2\delta\lambda)\omega_d(\kappa(x) + 1)K^d/\lambda,$$

which is less than $3\delta\omega_d(\kappa(x) + 1)K^d$ for sufficiently large $k$. Hence, $P[E^c] \to 0$ as $k \to \infty$. Moreover, by (3.5), we also have $P[F^c] \to 0$ as $k \to \infty$.

Assuming $\lambda$ to be so large that $|x - y| > 2\lambda^{-1/d}K$, if $E \cap F$ occurs, then

$$\mathcal{H}'_{\kappa(x)} \cap B^*_K = (\lambda^{1/d}(-x + \mathcal{P})) \cap B^*_K$$
$$= \lambda^{1/d}((-x + \mathcal{P}) \cap B^*_{\lambda^{-1/d}K}) = \lambda^{1/d}(-x + (\mathcal{P} \cap B^*_{\lambda^{-1/d}K}(x)))$$
$$= \lambda^{1/d}(-x + (\mathcal{X}'_m \cap B^*_{\lambda^{-1/d}K}(x)))$$
$$= \lambda^{1/d}((-x + \mathcal{X}'_m) \cap B^*_{\lambda^{-1/d}K}) = (\lambda^{1/d}(-x + \mathcal{X}'_m)) \cap B^*_K$$

so that $D(\mathcal{H}'_{\kappa(x)}, \lambda^{1/d}(-x + \mathcal{X}'_m)) \leq 1/K$. Hence, for any $K$, we have

$$P[D(\mathcal{H}'_{\kappa(x)}, \lambda^{1/d}(-x + \mathcal{X}'_m)) > 1/K] \to 0.$$

Similarly, we have

$$\max\{P[D(\mathcal{H}'_{\kappa(x)}, \lambda^{1/d}(-x + \mathcal{X}'_\ell)) > 1/K], P[D(\tilde{\mathcal{H}}'_{\kappa(y)}, \lambda^{1/d}(-y + \mathcal{X}'_m)) > 1/K],$$
$$P[D(\mathcal{H}'_{\kappa(x)}, \lambda^{1/d}(-x + (\mathcal{X}'_\ell)^{\mathbf{y}})) > 1/K], P[D(\mathcal{H}'_{\kappa(x)}, \lambda^{1/d}(-x + (\mathcal{X}'_m)^{\mathbf{y}})) > 1/K],$$
$$P[D(\tilde{\mathcal{H}}'_{\kappa(y)}, \lambda^{1/d}(-y + (\mathcal{X}'_m)^{\mathbf{x}})) > 1/K]\} \to 0.$$

Combining these, we have the required convergence in distribution.  □

**Lemma 3.3.** *Suppose $(x, y) \in \Omega_0 \times \Omega_0$, with $x \neq y$. Suppose, also, that $R(\mathbf{x}; \mathcal{H}_{\kappa(x)})$ and $R(\mathbf{y}; \mathcal{H}_{\kappa(y)})$ are almost surely finite. Suppose $(\lambda(m))_{m \geq 1}$ is a $(0, \infty)$-valued sequence with $\lambda(m)/m \to 1$ as $m \to \infty$. Then, for Borel $A \subseteq \mathbb{R}^d$, as $m \to \infty$, we have*

$$\xi_{\lambda(m)}(\mathbf{x}; \mathcal{X}_m, x + \lambda(m)^{-1/d} A) \xrightarrow{\mathcal{D}} \xi^{\mathbf{x}}_\infty(\mathcal{H}_{\kappa(x)}, A) \tag{3.6}$$

*and*

$$(\xi_{\lambda(m)}(\mathbf{x}; \mathcal{X}^{\mathbf{y}}_m, x + \lambda(m)^{-1/d} A), \xi_{\lambda(m)}(\mathbf{y}; \mathcal{X}^{\mathbf{x}}_m, y + \lambda(m)^{-1/d} A))$$
$$\xrightarrow{\mathcal{D}} (\xi^{\mathbf{x}}_\infty(\mathcal{H}_{\kappa(x)}, A), \xi^{\mathbf{y}}_\infty(\tilde{\mathcal{H}}_{\kappa(y)}, A)). \tag{3.7}$$

**Proof.** Given $A$, define the mapping $h_{A,x}: \mathcal{M} \times \mathcal{L} \to [0, \infty]$ and the mapping $h_A^2: (\mathcal{M} \times \mathcal{L} \times \mathcal{M} \times \mathcal{L}) \to [0, \infty]^2$ by

$$h_{A,x}(t, \mathcal{X}) = \xi((x, t); x + \mathcal{X}, x + A),$$
$$h_A^2(t, \mathcal{X}, t', \mathcal{X}') = (h_{A,x}(t, \mathcal{X}), h_{A,y}(t', \mathcal{X}')).$$

Since $R(\mathbf{x}; \mathcal{H}_{\kappa(x)}) < \infty$ a.s., the pair $(T, \mathcal{H}_{\kappa(x)})$ lies a.s. at a continuity point of $h_{A,x}$, where the topology on $\mathcal{M} \times \mathcal{L}$ is the product of the discrete topology on $\mathcal{M}$ and the topology induced by our metric $D$ on $\mathcal{L}$, defined at (3.1). Similarly, $(T, \mathcal{H}_{\kappa(x)}, T', \tilde{\mathcal{H}}_{\kappa(y)})$ lies a.s. at a continuity point of $h_A^2$. We have, by the definition of $\xi_\lambda$, that

$$\xi_\lambda(\mathbf{x}; \mathcal{X}_m, x + \lambda^{-1/d} A)$$
$$= h_{A,x}(T, \lambda^{1/d}(-x + \mathcal{X}_m)),$$
$$(\xi_\lambda(\mathbf{x}; \mathcal{X}^{\mathbf{y}}_m, x + \lambda^{-1/d} A), \xi_\lambda(\mathbf{y}; \mathcal{X}^{\mathbf{x}}_m, y + \lambda^{-1/d} A))$$
$$= h_A^2(T, \lambda^{1/d}(-x + \mathcal{X}^{\mathbf{y}}_m), T', \lambda^{1/d}(-y + \mathcal{X}^{\mathbf{x}}_m)).$$

By Lemma 3.2, we have $(T, \lambda(m)^{1/d}(-x + \mathcal{X}_m)) \xrightarrow{\mathcal{D}} (T, \mathcal{H}_{\kappa(x)})$ so that (3.6) follows by the Continuous Mapping Theorem ([3], Chapter 1, Theorem 5.1). Also, by Lemma 3.2,

$$(T, \lambda(m)^{1/d}(-x + \mathcal{X}^{\mathbf{y}}_m), T', \lambda(m)^{1/d}(-y + \mathcal{X}^{\mathbf{x}}_m)) \xrightarrow{\mathcal{D}} (T, \mathcal{H}_{\kappa(x)}, T', \tilde{\mathcal{H}}_{\kappa(y)})$$

so that (3.7) also follows by the Continuous Mapping Theorem. □

**Lemma 3.4.** *Suppose* $(x,y) \in \Omega_0 \times \Omega_0$, *with* $x \neq y$, *that* $R(\mathbf{x}; \mathcal{H}_{\kappa(x)})$, $R(\mathbf{y}; \mathcal{H}_{\kappa(y)})$, $\xi_\infty^{\mathbf{x}}(\mathcal{H}_{\kappa(x)}, \mathbb{R}^d)$ *and* $\xi_\infty^{\mathbf{y}}(\mathcal{H}_{\kappa(y)}, \mathbb{R}^d)$ *are almost surely finite and that* $(\lambda(m))_{m \geq 1}$ *is a* $(0, \infty)$-*valued sequence with* $\lambda(m)/m \to 1$ *as* $m \to \infty$. *Then, as* $m \to \infty$, *we have*

$$\xi_{\lambda(m)}(\mathbf{x}; \mathcal{X}_m, \Omega) \xrightarrow{\mathcal{D}} \xi_\infty^{\mathbf{x}}(\mathcal{H}_{\kappa(x)}, \mathbb{R}^d) \tag{3.8}$$

*and*

$$(\xi_{\lambda(m)}(\mathbf{x}; \mathcal{X}_m^{\mathbf{y}}, \Omega), \xi_{\lambda(m)}(\mathbf{y}; \mathcal{X}_m^{\mathbf{x}}, \Omega)) \xrightarrow{\mathcal{D}} (\xi_\infty^{\mathbf{x}}(\mathcal{H}_{\kappa(x)}, \mathbb{R}^d), \xi_\infty^{\mathbf{y}}(\tilde{\mathcal{H}}_{\kappa(y)}, \mathbb{R}^d)). \tag{3.9}$$

**Proof.** Since $\Omega_0 \subseteq \Omega_\infty \subseteq \Omega$ and $\Omega_\infty$ is open, for any $K > 0$, we have, for sufficiently large $m$, that

$$0 \leq \xi_{\lambda(m)}(\mathbf{x}; \mathcal{X}_m, \mathbb{R}^d \setminus \Omega)$$
$$\leq \xi_{\lambda(m)}(\mathbf{x}; \mathcal{X}_m, \mathbb{R}^d \setminus B_{\lambda(m)^{-1/d}K}(x)) \xrightarrow{\mathcal{D}} \xi_\infty^{\mathbf{x}}(\mathcal{H}_{\kappa(x)}, \mathbb{R}^d \setminus B_K), \tag{3.10}$$

where the convergence follows from (3.6). By assumption, $\xi_\infty^{\mathbf{x},T}(\mathcal{H}_{\kappa(\mathbf{x})}, \mathbb{R}^d)$ is almost surely finite, so the limit in (3.10) itself tends to zero in probability as $K \to \infty$ and therefore $\xi_{\lambda(m)}(\mathbf{x}; \mathcal{X}_m, \mathbb{R}^d \setminus \Omega) \xrightarrow{P} 0$ as $\lambda \to \infty$. Combining this with the case $A = \mathbb{R}^d$ of (3.6) and using Slutsky's theorem (see, e.g., [13]), we obtain (3.8).

A similar argument to the above, using (3.7), shows that as $m \to \infty$,

$$(\xi_{\lambda(m)}(\mathbf{x}; \mathcal{X}_m^{\mathbf{y}}, \mathbb{R}^d \setminus \Omega), \xi_{\lambda(m)}(\mathbf{y}; \mathcal{X}_m^{\mathbf{x}}, \mathbb{R}^d \setminus \Omega)) \xrightarrow{P} (0, 0)$$

and by using this with the case $A = \mathbb{R}^d$ of (3.7) and Slutsky's theorem in two dimensions, we obtain (3.9). □

The next lemma compares the measure $\xi_\lambda(\mathbf{x}; \mathcal{X}, \cdot)$ to the corresponding point measure $\xi_\lambda^*(\mathbf{x}; \mathcal{X}, \cdot)$. In proving this, for $f \in B(\Omega)$, we write $\|f\|_\infty$ for $\sup\{|f(x)| : x \in \Omega\}$.

**Lemma 3.5.** *Let* $x \in \Omega_0$ *and suppose that* $R(\mathbf{x}; \mathcal{H}_{\kappa(x)})$ *and* $\xi_\infty^{\mathbf{x}}(\mathcal{H}_{\kappa(x)}, \mathbb{R}^d)$ *are almost surely finite. Let* $y \in \mathbb{R}^d$, *with* $y \neq x$. *Suppose that* $f \in B(\Omega)$ *and suppose either that* $f$ *is continuous at* $x$ *or that assumption* A2 *holds and* $x$ *is a Lebesgue point of* $f$. *Suppose* $(\lambda(m))_{m \geq 1}$ *is a* $(0, \infty)$-*valued sequence with* $\lambda(m)/m \to 1$ *as* $m \to \infty$. *Then, as* $m \to \infty$,

$$\langle f, \xi_{\lambda(m)}(\mathbf{x}; \mathcal{X}_m) - \xi_{\lambda(m)}^*(\mathbf{x}; \mathcal{X}_m) \rangle \xrightarrow{P} 0 \tag{3.11}$$

*and*

$$\langle f, \xi_{\lambda(m)}(\mathbf{x}; \mathcal{X}_m^{\mathbf{y}}) - \xi_{\lambda(m)}^*(\mathbf{x}; \mathcal{X}_m^{\mathbf{y}}) \rangle \xrightarrow{P} 0. \tag{3.12}$$

**Proof.** In this proof, we write $\lambda$ for $\lambda(m)$. The left-hand side of (3.11) is equal to

$$\int_{\mathbb{R}^d} (f(z) - f(x))\xi_\lambda(\mathbf{x}; \mathcal{X}_m, \mathrm{d}z). \tag{3.13}$$

Given $K > 0$, we split the region of integration in (3.13) into the complementary regions $B_{\lambda^{-1/d}K}(x)$ and $\mathbb{R}^d \setminus B_{\lambda^{-1/d}K}(x)$. Consider the latter region first. By (3.6), we have

$$\left|\int_{\mathbb{R}^d \setminus B_{\lambda^{-1/d}K}(x)} (f(z) - f(x))\xi_\lambda(\mathbf{x}; \mathcal{X}_m, \mathrm{d}z)\right| \leq 2\|f\|_\infty \xi_\lambda(\mathbf{x}; \mathcal{X}_m, \mathbb{R}^d \setminus B_{\lambda^{-1/d}K}(x))$$

$$\xrightarrow{\mathcal{D}} 2\|f\|_\infty \xi_\infty^{\mathbf{x}}(\mathcal{H}_{\kappa(x)}, \mathbb{R}^d \setminus B_K),$$

where the limit is almost surely finite and converges in probability to zero as $K \to \infty$. Hence for $\varepsilon > 0$, we have

$$\lim_{K \to \infty} \limsup_{m \to \infty} P\left[\left|\int_{\mathbb{R}^d \setminus B_{\lambda^{-1/d}K}(x)} (f(z) - f(x))\xi_\lambda(\mathbf{x}; \mathcal{X}_m, \mathrm{d}z)\right| > \varepsilon\right] = 0. \tag{3.14}$$

Turning to the integral over $B_{\lambda^{-1/d}K}(x)$, we consider separately the case where $f$ is continuous at $x$ and the case where A2 holds and $x$ is a Lebesgue point of $f$. To deal with the first of these cases, writing $\phi_\varepsilon(x)$ for $\sup\{|f(y) - f(x)| : y \in B_\varepsilon(x)\}$, we observe that

$$\left|\int_{B_{\lambda^{-1/d}K}(x)} (f(z) - f(x))\xi_\lambda(\mathbf{x}; \mathcal{X}_m, \mathrm{d}z)\right| \leq \phi_{\lambda^{-1/d}K}(x)\xi_\lambda(\mathbf{x}; \mathcal{X}_m, \Omega). \tag{3.15}$$

If $f$ is continuous at $x$, then $\phi_{\lambda^{-1/d}K}(x) \to 0$, while $\xi_\lambda(\mathbf{x}; \mathcal{X}_m, \Omega)$ converges in distribution to the finite random variable $\xi_\infty^{\mathbf{x}}(\mathcal{H}_{\kappa(x)}, \mathbb{R}^d)$, by (3.8), and hence the right-hand side of (3.15) tends to zero in probability as $m \to \infty$. Combined with (3.14), this gives us (3.11) in the case where $f$ is continuous at $x$.

Under assumption A2, for Borel $A \subseteq \mathbb{R}^d$, the change of variables $z = x + \lambda^{-1/d}(y - x)$ yields

$$\xi_\lambda(\mathbf{x}; \mathcal{X}, A) = \int_{x + \lambda^{1/d}(-x+A)} \xi'(\mathbf{x}; x + \lambda^{1/d}(-x+\mathcal{X}), y)\,\mathrm{d}y$$

$$= \lambda \int_A \xi'(\mathbf{x}; x + \lambda^{1/d}(-x+\mathcal{X}), x + \lambda^{1/d}(z-x))\,\mathrm{d}z.$$

Hence, under A2,

$$\left|\int_{B_{\lambda^{-1/d}K}(x)} (f(z) - f(x))\xi_\lambda(x; \mathcal{X}_m, \mathrm{d}z)\right|$$

$$= \lambda\left|\int_{B_{\lambda^{-1/d}K}(x)} (f(z) - f(x))\xi'(\mathbf{x}; x + \lambda^{1/d}(-x+\mathcal{X}_m), x + \lambda^{1/d}(z-x))\,\mathrm{d}z\right|$$

$$\leq K_0 \lambda \int_{B_{\lambda^{-1/d}K}(x)} |f(z) - f(x)| \, \mathrm{d}z$$

and if, additionally, $x$ is a Lebesgue point of $f$, then this tends to zero. Combined with (3.14), this gives us (3.11) in the case where A2 holds and $x$ is a Lebesgue point of $f$.

The proof of (3.12) is similar; we use (3.7) and (3.9) instead of (3.6) and (3.8). □

By combining Lemmas 3.3 and 3.5, we obtain the following, which is the main ingredient in our proof of the Law of Large Numbers in Theorem 2.1.

**Lemma 3.6.** *Suppose $(x,y) \in \Omega_0 \times \Omega_0$, with $x \neq y$. Suppose, also, that $R(\mathbf{x}; \mathcal{H}_{\kappa(x)})$, $R(\mathbf{y}; \mathcal{H}_{\kappa(y)})$, $\xi_\infty^{\mathbf{x}}(\mathcal{H}_{\kappa(x)}, \mathbb{R}^d)$ and $\xi_\infty^{\mathbf{y}}(\mathcal{H}_{\kappa(y)}, \mathbb{R}^d)$ are almost surely finite. Let $f \in B(\Omega)$ and suppose either that A1 holds, that A2 holds and $x$ is a Lebesgue point of $f$, or that $f$ is continuous at $x$. Suppose $(\lambda(m))_{m \geq 1}$ is a $(0, \infty)$-valued sequence with $\lambda(m)/m \to 1$ as $m \to \infty$. Then, as $m \to \infty$,*

$$\langle f, \xi_{\lambda(m)}(\mathbf{x}; \mathcal{X}_m) \rangle \xrightarrow{\mathcal{D}} f(x) \xi_\infty^{\mathbf{x}}(\mathcal{H}_{\kappa(x)}, \mathbb{R}^d) \qquad (3.16)$$

*and*

$$\langle f, \xi_{\lambda(m)}(\mathbf{x}; \mathcal{X}_m^{\mathbf{y}}) \rangle \langle f, \xi_{\lambda(m)}(\mathbf{y}; \mathcal{X}_m^{\mathbf{x}}) \rangle$$
$$\xrightarrow{\mathcal{D}} f(x) f(y) \xi_\infty^{\mathbf{x}}(\mathcal{H}_{\kappa(x)}, \mathbb{R}^d) \xi_\infty^{\mathbf{y}}(\tilde{\mathcal{H}}_{\kappa(y)}, \mathbb{R}^d). \qquad (3.17)$$

**Proof.** Note, first, that by (3.8),

$$\langle f, \xi_{\lambda(m)}^*(\mathbf{x}; \mathcal{X}_m) \rangle = f(x) \xi_{\lambda(m)}(\mathbf{x}; \mathcal{X}_m, \Omega) \xrightarrow{\mathcal{D}} f(x) \xi_\infty^{\mathbf{x}}(\mathcal{H}_{\kappa(x)}, \mathbb{R}^d) \qquad (3.18)$$

and similarly, by (3.9),

$$(\langle f, \xi_{\lambda(m)}^*(\mathbf{x}; \mathcal{X}_m^{\mathbf{y}}) \rangle, \langle f, \xi_{\lambda(m)}^*(\mathbf{y}; \mathcal{X}_m^{\mathbf{x}}) \rangle)$$
$$\xrightarrow{\mathcal{D}} (f(x) \xi_\infty^{\mathbf{x}}(\mathcal{H}_{\kappa(x)}, \mathbb{R}^d), f(y) \xi_\infty^{\mathbf{y}}(\tilde{\mathcal{H}}_{\kappa(y)}, \mathbb{R}^d)). \qquad (3.19)$$

In the case where A1 holds, we have $\xi_\lambda = \xi_\lambda^*$, so (3.16) follows immediately from (3.18) and (3.17) follows immediately from (3.19).

In the other two cases described, we have (3.11), by Lemma 3.5. Combining this with (3.18), we see, by Slutsky's theorem, that (3.16) still holds in the other two cases. Similarly, by (3.19), (3.12) and Slutsky's theorem, we can obtain (3.17) in the other cases too. □

## 4. Proof of general laws of large numbers

In this section, we complete the proofs of Theorems 2.1 and 2.2, using the weak convergence results from the preceding section. Throughout this section, we assume that $(\lambda(n))_{n\geq 1}$ satisfy $\lambda(n) > 0$ and that $\lambda(n)/n \to 1$ as $n \to \infty$. Also, let $\mathcal{H}_{\kappa(X)}$ denote a Cox point process in $\mathbb{R}^d \times \mathcal{M}$, whose distribution, given $X = x$, is that of $\mathcal{H}_{\kappa(x)}$ (where $\mathbf{X} = (X, T)$ is as in Section 2). We first show that the conditions of Theorem 2.1 imply that $\xi_\infty^{(x,T)}(\mathcal{H}_{\kappa(x)}, \mathbb{R}^d)$ is finite.

**Lemma 4.1.** *Suppose* $\mathbb{R}((x,T); \mathcal{H}_{\kappa(x)}) < \infty$ *almost surely, for $\kappa$-almost all $x \in \Omega_\infty$. If (2.5) holds for $p = 1$, then $\xi_\infty^{(x,T)}(\mathcal{H}_{\kappa(x)}, \mathbb{R}^d)$ is almost surely finite, for $\kappa$-almost all $x \in \Omega_\infty$.*

**Proof.** Given $K > 0$, define the random variables

$$S_K := \xi_\infty^{\mathbf{X}}(\mathcal{H}_{\kappa(X)}, B_K)\mathbf{1}_{\Omega_\infty}(X), \qquad S_{K,n} := \xi_{\lambda(n)}(\mathbf{X}; \mathcal{X}_{n-1}, B_{\lambda(n)^{-1/d}K}(X)).$$

By Lemma 3.3, for any bounded continuous test function $h$ on $\mathbb{R}$, as $n \to \infty$, we have almost sure convergence of $\mathbb{E}[h(S_{K,n})|\mathbf{X}]$ to $\mathbb{E}[h(S_K)|\mathbf{X}]$. By taking expectations and using the dominated convergence theorem, we have that $\mathbb{E}[h(S_{K,n})] \to \mathbb{E}[h(S_K)]$. Hence, $S_{K,n} \xrightarrow{\mathcal{D}} S_K$ as $n \to \infty$. Hence, by (2.5) and Fatou's Lemma, $\mathbb{E}[S_K]$ is bounded by a constant independent of $K$. Taking $K \to \infty$, we may deduce that $\xi_\infty^{\mathbf{X}}(\mathcal{H}_{\kappa(X)}, \mathbb{R}^d)\mathbf{1}_{\Omega_\infty}(X)$ has finite mean and so is almost surely finite. The result follows. $\square$

To prove Theorem 2.1, we shall use the following general expressions for the first two moments of $\langle f, \nu_{\lambda,n} \rangle$. By (2.3), we have

$$n^{-1}\mathbb{E}\langle f, \nu_{\lambda,n}\rangle = \mathbb{E}\langle f, \xi_\lambda(\mathbf{X}; \mathcal{X}_{n-1})\rangle \tag{4.1}$$

and

$$n^{-2} \operatorname{Var}\langle f, \nu_{\lambda(n),n}^\xi\rangle = n^{-1}\mathbb{E}[\langle f, \xi_{\lambda(n)}(\mathbf{X}; \mathcal{X}_{n-1})\rangle^2]$$
$$+ \left(\frac{n-1}{n}\right)\mathbb{E}[\langle f, \xi_{\lambda(n)}(\mathbf{X}; \mathcal{X}_{n-2}^{\mathbf{X}'})\rangle\langle f, \xi_{\lambda(n)}(\mathbf{X}'; \mathcal{X}_{n-2}^{\mathbf{X}})\rangle]$$
$$- (\mathbb{E}[\langle f, \xi_{\lambda(n)}(\mathbf{X}; \mathcal{X}_{n-1})\rangle])^2. \tag{4.2}$$

Recall that by definition (2.1), $\xi_\lambda((x,t); \mathcal{X}, \mathbb{R}^d) = 0$ for $x \in \mathbb{R}^d \setminus \Omega_\lambda$, with $(\Omega_\lambda, \lambda \geq 1)$ a given non-decreasing family of Borel subsets of $\mathbb{R}^d$ with limit set $\Omega_\infty$ and $\Omega_\infty \subseteq \Omega \subseteq \mathbb{R}^d$. In the simplest case, $\Omega_\lambda = \mathbb{R}^d$ for all $\lambda$.

**Proof of Theorem 2.1.** Let $f \in B(\Omega)$. First, we prove (i) for the case $q = 2$. Assume that (2.5) holds for some $p > 2$. Set $J := f(X)\xi_\infty^{\mathbf{X}}(\mathcal{H}_{\kappa(X)}, \mathbb{R}^d)\mathbf{1}_{\Omega_\infty}(X)$ and let $J'$ be an independent copy of $J$. By Lemma 4.1, $J$ is almost surely finite.

For any bounded continuous test function $h$ on $\mathbb{R}$, by (3.16) from Lemma 3.6, as $n \to \infty$, we have $\mathbb{E}[h(\langle f, \xi_{\lambda(n)}(\mathbf{X}; \mathcal{X}_{n-1})\rangle)|\mathbf{X}] \to \mathbb{E}[J|\mathbf{X}]$, almost surely. Hence, $\mathbb{E}[h(\langle f, \xi_{\lambda(n)}(\mathbf{X}; \mathcal{X}_{n-1})\rangle)] \to \mathbb{E}[h(J)]$ so that

$$\langle f, \xi_{\lambda(n)}(\mathbf{X}; \mathcal{X}_{n-1})\rangle \xrightarrow{\mathcal{D}} J. \tag{4.3}$$

*LLNs in stochastic geometry with applications*　　　　　　　　　　　　　　　　　　1139

Similarly, using (3.17), we obtain

$$\langle f, \xi_{\lambda(n)}(\mathbf{X}; \mathcal{X}_{n-2}^{\mathbf{X}'})\rangle \langle f, \xi_{\lambda(n)}(\mathbf{X}'; \mathcal{X}_{n-2}^{\mathbf{X}})\rangle \xrightarrow{\mathcal{D}} J'J. \qquad (4.4)$$

Also, by (2.5) and the Cauchy–Schwarz inequality, the variables in the left-hand side of (4.3) and in the left-hand side of (4.4) are uniformly integrable, so we have convergence of means in both cases. Also, (2.5) shows that the first term in the right-hand side of (4.2) tends to zero. Hence, we find that the expression (4.2) tends to zero. Moreover, by (4.1) and the convergence of expectations corresponding to (4.3), $n^{-1}\mathbb{E}\langle f, \nu_{\lambda(n),n}\rangle$ tends to $\mathbb{E}[J]$ and this gives us (2.6) with $L^2$ convergence,

Now, consider the case $q = 1$. Assume (2.5) holds for some $p > 1$. First, assume $f$ is non-negative. We use a truncation argument; for $K > 0$, let $\xi_\lambda^K$ be the truncated version of the measure $\xi_\lambda^*$, defined by

$$\xi_\lambda^K((x,t); \mathcal{X}, A) := \min(\xi_\lambda((x,t); \mathcal{X}, \Omega), K)\mathbf{1}_A(x).$$

Let $\Omega^*$ be the set of $x \in \Omega_0$ such that $R((x,T); \mathcal{H}_{\kappa(x)})$ and $\xi_\infty(\mathcal{H}_{\kappa(x)}, \mathbb{R}^d)$ are almost surely finite. By Lemma 4.1, $\kappa(\Omega_0 \setminus \Omega^*) = 0$.

Then, for $\mathbf{x} = (x, T)$ with $x \in \Omega^*$,

$$\langle f, \xi_{\lambda(n)}^K(\mathbf{x}; \mathcal{X}_{n-1})\rangle = f(x)\min(\xi_\lambda(\mathbf{x}; \mathcal{X}_{n-1}, \Omega), K)$$
$$\xrightarrow{\mathcal{D}} f(x)\min(\xi_\infty^{\mathbf{x}}(\mathcal{H}_{\kappa(x)}, \mathbb{R}^d), K), \qquad (4.5)$$

where the convergence follows from (3.8). Similarly, for distinct $x, y$ in $\Omega^*$, setting $\mathbf{x} = (x,T)$ and $\mathbf{y} = (y, T')$, by (3.9), we have that

$$\langle f, \xi_{\lambda(n)}^K(\mathbf{x}; \mathcal{X}_{n-2}^{\mathbf{y}})\rangle \langle f, \xi_{\lambda(n)}^K(\mathbf{y}; \mathcal{X}_{n-2}^{\mathbf{x}})\rangle$$
$$\xrightarrow{\mathcal{D}} f(x)f(y)\min(\xi_\infty^{\mathbf{x}}(\mathcal{H}_{\kappa(x)}, \mathbb{R}^d), K)\min(\xi_\infty^{\mathbf{y}}(\tilde{\mathcal{H}}_{\kappa(y)}, \mathbb{R}^d), K). \qquad (4.6)$$

Using (4.5) and the same argument as for (4.3), we may deduce that

$$\langle f, \xi_{\lambda(n)}^K(\mathbf{X}; \mathcal{X}_{n-1})\rangle \xrightarrow{\mathcal{D}} J_K, \qquad (4.7)$$

where we set

$$J_K := f(X)\min(\xi_\infty^{\mathbf{X}}(\mathcal{H}_{\kappa(X)}, \mathbb{R}^d), K)\mathbf{1}_\Omega(X).$$

Likewise, using (4.6), we obtain

$$\langle f, \xi_\lambda^K(\mathbf{X}; \mathcal{X}_{n-2}^{\mathbf{X}'})\rangle \langle f, \xi_\lambda^K(\mathbf{X}'; \mathcal{X}_{n-2}^{\mathbf{X}})\rangle \xrightarrow{\mathcal{D}} J_K J_K', \qquad (4.8)$$

where $J_K'$ is an independent copy of $J_K$. Also, since $\xi_\lambda^K(\mathbf{x}; \mathcal{X}, \Omega)$ is bounded by $K$, the distributional convergences (4.7) and (4.8) are of bounded variables, so the corresponding convergence of expectations holds. Set

$$\nu^K_{\lambda,n} := \sum_{i=1}^n \xi^K_\lambda(\mathbf{X}_i; \mathcal{X}_n), \qquad \nu^*_{\lambda,n} := \sum_{i=1}^n \xi^*_\lambda(\mathbf{X}_i; \mathcal{X}_n).$$

By following the proof of (2.6) with $L^2$ convergence, we obtain

$$n^{-1}\langle f, \nu^K_{\lambda(n),n}\rangle \xrightarrow{L^2} \mathbb{E} J_K. \tag{4.9}$$

Also,

$$0 \leq \mathbb{E}[n^{-1}\langle f, \nu^*_{\lambda(n),n}\rangle - n^{-1}\langle f, \nu^K_{\lambda(n),n}\rangle]$$
$$= \mathbb{E}[\langle f, \xi^*_{\lambda(n)}(\mathbf{X}; \mathcal{X}_{n-1}) - \xi^K_{\lambda(n)}(\mathbf{X}; \mathcal{X}_{n-1})\rangle]$$
$$\leq \|f\|_\infty \mathbb{E}[\xi_{\lambda(n)}(\mathbf{X}; \mathcal{X}_{n-1}, \Omega)\mathbf{1}\{\xi_{\lambda(n)}(\mathbf{X}; \mathcal{X}_{n-1}, \Omega) > K\}],$$

which tends to zero as $K \to \infty$, uniformly in $n$, because the moments condition (2.5), $p > 1$, implies that the random variables $\xi_{\lambda(n)}(\mathbf{X}; \mathcal{X}_{n-1}, \Omega)$ are uniformly integrable. Also, by monotone convergence, as $K \to \infty$, the right-hand side of (4.9) converges to $\mathbb{E}[J]$. Hence, taking $K \to \infty$ in (4.9) yields

$$n^{-1}\langle f, \nu^*_{\lambda(n),n}\rangle \xrightarrow{L^1} \mathbb{E}[J]. \tag{4.10}$$

This gives us (2.6) with $L^1$ convergence when assumption A1 holds, in the case where $f$ is non-negative; by taking positive and negative parts of $f$ and using linearity, we can extend this to general $f$.

Now, suppose A2 or A3 holds. Then,

$$\mathbb{E}\left[n^{-1}\sum_{i=1}^n |\langle f, \xi_{\lambda(n)}(\mathbf{X}_i; \mathcal{X}_n) - \xi^*_{\lambda(n)}(\mathbf{X}_i; \mathcal{X}_n)\rangle|\right]$$
$$= \mathbb{E}[|\langle f, \xi_{\lambda(n)}(\mathbf{X}; \mathcal{X}_{n-1}) - \xi^*_{\lambda(n)}(\mathbf{X}; \mathcal{X}_{n-1})\rangle|]. \tag{4.11}$$

By (3.11), the variables $|\langle f, \xi_{\lambda(n)}(\mathbf{X}; \mathcal{X}_{n-1}) - \xi^*_{\lambda(n)}(\mathbf{X}; \mathcal{X}_{n-1})\rangle|$ tend to zero in probability and, by (2.5), they are uniformly integrable, so their mean tends to zero, that is, the expression (4.11) tends to zero and thus we have (2.7). Combining this with (4.10) gives us (2.6) for $q = 1$ when assumption A2 or A3 holds, completing the proof. □

**Proof of Theorem 2.2.** Suppose that $R((x,T); \mathcal{H}_{\kappa(x)})$ is almost surely finite for $\kappa$-almost all $x \in \Omega_\infty$. Suppose that $\lambda(n)/n \to 1$ as $n \to \infty$ and that there exists $p > 1$ such that (2.5) holds. By the case $p = 1$ of Theorem 2.1 (or, more directly, by the argument at the start of the proof of that result), we have convergence of means in (2.6). To derive almost sure convergence under condition (2.9), we loosely follow the argument from [17], pages 298–299. For $\lambda > 0$, define $H_\lambda : \bigcup_{n=1}^{\infty} [(\mathbb{R}^d \times \mathcal{M})^n] \to \mathbb{R}$ by

$$H_\lambda(\mathbf{x}_1, \ldots, \mathbf{x}_n) := \sum_{i=1}^{n} \langle f, \xi_\lambda(\mathbf{x}_i; \{\mathbf{x}_1, \ldots, \mathbf{x}_n\}) \rangle.$$

Then, $\langle f, \nu_{\lambda,n} \rangle = H_\lambda(\mathbf{X}_1, \ldots, \mathbf{X}_n)$. Let $\mathcal{F}_i$ denote the $\sigma$-field generated by $\mathbf{X}_1, \ldots, \mathbf{X}_i$ and let $\mathcal{F}_0$ denote the trivial $\sigma$-field. We then have the martingale difference representation $\langle f, \nu_{\lambda(n),n} \rangle - \mathbb{E}\langle f, \nu_{\lambda(n),n} \rangle = \sum_{i=1}^{n} d_i$, where $d_i := \mathbb{E}[\langle f, \nu_{\lambda(n),n} \rangle | \mathcal{F}_i] - \mathbb{E}[\langle f, \nu_{\lambda(n),n} \rangle | \mathcal{F}_{i-1}]$. Notice that

$$d_i = \mathbb{E}[H_{\lambda(n)}(\mathbf{X}_1, \ldots, \mathbf{X}_n) - H_{\lambda(n)}(\mathbf{X}_1, \ldots, \mathbf{X}_{i-1}, \mathbf{X}', \mathbf{X}_{i+1}, \ldots, \mathbf{X}_n) | \mathcal{F}_i].$$

By assumption (2.9), $\|H_{\lambda(n)}(\mathbf{X}_1, \ldots, \mathbf{X}_n) - H_{\lambda(n)}(\mathbf{X}_1, \ldots, \mathbf{X}_{n-1})\|_{p'}$ is bounded by a constant $C$, independent of $n$, and so, by Minkowski's inequality and exchangeability of $\mathbf{X}_1, \ldots, \mathbf{X}_n$,

$$\|H_{\lambda(n)}(\mathbf{X}_1, \ldots, \mathbf{X}_n) - H_{\lambda(n)}(\mathbf{X}_1, \ldots, \mathbf{X}_{i-1}, \mathbf{X}', \mathbf{X}_{i+1}, \ldots, \mathbf{X}_n)\|_{p'}$$
$$\leq 2\|H_{\lambda(n)}(\mathbf{X}_1, \ldots, \mathbf{X}_n) - H_{\lambda(n)}(\mathbf{X}_1, \ldots, \mathbf{X}_{n-1})\|_{p'} \leq 2C,$$

so, by the conditional Jensen inequality, allowing the constant $C$ to change from line to line, we have

$$\mathbb{E}|d_i|^{p'} \leq \mathbb{E}\mathbb{E}[|H_{\lambda(n)}(\mathbf{X}_1, \ldots, \mathbf{X}_n) - H_{\lambda(n)}(\mathbf{X}_1, \ldots, \mathbf{X}_{i-1}, \mathbf{X}', \mathbf{X}_{i+1}, \ldots, \mathbf{X}_n)|^{p'} | \mathcal{F}_i]$$
$$\leq C. \tag{4.12}$$

Choose $\gamma$ to satisfy $\gamma < 1/2$ and $p'\gamma > \beta + 1$. By the condition $p' > 2(\beta + 1)$, such $\gamma$ exists.

We now use the following modification of Azuma's inequality, introduced by Chalker *et al.* ([4], Lemma 1). For any martingale difference sequence $d_i$, $i \geq 1$, and for all sequences $w_i$, $i \geq 1$, of positive numbers, we have, for all $t > 0$, that

$$P\left[\left|\sum_{i=1}^{n} d_i\right| > t\right] \leq 2\exp\left(\frac{-t^2}{32 \sum_{i=1}^{n} w_i^2}\right)$$
$$+ \left(1 + 2t^{-1} \sup_i \|d_i\|_\infty\right) \sum_{i=1}^{n} P[|d_i| > w_i].$$

Letting $w_i := n^\gamma$, $t := \varepsilon n$, using (4.12) and Markov's inequality and noting that $\sup_i \|d_i\|_\infty \leq Cn^\beta$ by the first part of (2.9), we obtain, for any $\varepsilon > 0$, that

$$P\left[\left|\sum_{i=1}^{n} d_i\right| > \varepsilon n\right] \leq 2\exp\left(\frac{-n^2}{Cn^{1+2\gamma}}\right) + (1 + Cn^{\beta-1})\frac{n}{n^{p'\gamma}},$$

which is summable in $n$ by the choice of $\gamma$ (since we assume $\beta \geq 1$). Hence, by the Borel–Cantelli lemma, we have almost sure convergence for (2.6).

To prove (2.7) with almost sure convergence under assumption (2.10), define $\tilde{H}_\lambda : \bigcup_{n=1}^{\infty} [(\mathbb{R}^d \times \mathcal{M})^n] \to \mathbb{R}$ by

$$\tilde{H}_\lambda(\mathbf{x}_1, \ldots, \mathbf{x}_n) := \sum_{i=1}^{n} |\langle f, \xi_{\lambda(n)}(\mathbf{x}_i; \{\mathbf{x}_1, \ldots, \mathbf{x}_n\}) - \xi^*_{\lambda(n)}(\mathbf{x}_i; \{\mathbf{x}_1, \ldots, \mathbf{x}_n\}) \rangle|$$

and then follow the same argument as given above, with $H_\lambda$ replaced by $\tilde{H}_\lambda$.

Next, we show that (2.11) implies (2.9) and (2.10) for any $f \in B(\Omega)$. Since

$$\Delta_{\lambda,n}(f) = \langle f, \xi_\lambda(\mathbf{X}_n; \mathcal{X}_n) \rangle + \sum_{i=1}^{n-1} \langle f, \xi_\lambda(\mathbf{X}_i; \mathcal{X}_n) - \xi_\lambda(\mathbf{X}_i; \mathcal{X}_{n-1}) \rangle$$

and since, for any signed measure $\mu$ on $\Omega$ with total variation $|\mu|$, we have $\langle f, \mu \rangle \le \|f\|_\infty \times |\mu|$, it follows by the triangle inequality and definition (2.8) that $|\Delta_{\lambda,n}(f)| \le \|f\|_\infty \tilde{\Delta}_{\lambda,n}$ and hence (2.11) implies (2.9).

Finally, we show that (2.11) implies (2.10). By definition,

$$\begin{aligned} \Delta'_{\lambda,n}(f) &= G_{\lambda,n}(f) - G_{\lambda,n-1}(f) \\ &= |\langle f, \xi_\lambda(\mathbf{X}_n; \mathcal{X}_n) - \xi^*_\lambda(\mathbf{X}_n; \mathcal{X}_n) \rangle| \\ &\quad + \sum_{i=1}^{n-1} (|\langle f, \xi_\lambda(\mathbf{X}_i; \mathcal{X}_n) - \xi^*_\lambda(\mathbf{X}_i; \mathcal{X}_n) \rangle| \\ &\quad - |\langle f, \xi_\lambda(\mathbf{X}_i; \mathcal{X}_{n-1}) - \xi^*_\lambda(\mathbf{X}_i; \mathcal{X}_{n-1}) \rangle|). \end{aligned} \quad (4.13)$$

For any real $a_1, a_2, b_1, b_2$, we have $|(|a_1 - b_1| - |a_2 - b_2|)| \le |a_1 - a_2| + |b_1 - b_2|$, by the triangle inequality, and using this, we can deduce from (4.13) that

$$|\Delta'_{\lambda,n}(f)| \le 2\|f\|_\infty \xi_\lambda(\mathbf{X}_n; \mathcal{X}_n, \Omega) + 4\|f\|_\infty \sum_{i=1}^{n-1} |\xi_\lambda(\mathbf{X}_i; \mathcal{X}_n)|(\Omega).$$

By definition (2.8), this is at most $4\|f\|_\infty \tilde{\Delta}_{\lambda,n}(f)$, so (2.11) implies (2.10). □

## 5. Applications of the general theory

### 5.1. Voronoi estimation of a set

The first example illustrating our general result is concerned with coverage of a set by Voronoi cells. Let $\Omega := (0,1)^d$. For finite $\mathcal{X} \subset \mathbb{R}^d$ and $\mathbf{x} \in \mathcal{X}$, let $\tilde{V}(\mathbf{x}; \mathcal{X})$ denote the closed Voronoi cell with nucleus $\mathbf{x}$ for the Voronoi tessellation induced by $\mathcal{X}$, that is, the set of $\mathbf{y} \in \mathbb{R}^d$ lying at least as close to $\mathbf{x}$ (in the Euclidean sense) as to any other point of $\mathcal{X}$. Let $V(\mathbf{x}; \mathcal{X})$ be the intersection of $\tilde{V}(\mathbf{x}; \mathcal{X})$ with $\Omega$. Let $\kappa$ be a density function on $\Omega$, let $\mathbf{X}_1, \mathbf{X}_2, \ldots$ be independent random $d$-vectors taking values in $\Omega$ with common probability density $\kappa$ and let $\mathcal{X}_n = \{\mathbf{X}_1, \ldots, \mathbf{X}_n\}$ (in this section, boldface vectors represent unmarked points in $\mathbb{R}^d$).

Let $A$ be an arbitrary Borel subset of $\Omega$. Let $A_n$ be the estimator of the (unknown) set $A$ based on data from sensors at $\mathcal{X}_n$ using Voronoi cells, that is, let $A_n := \bigcup_{\mathbf{x} \in \mathcal{X}_n \cap A} V(\mathbf{x}; \mathcal{X}_n)$. With a view to potential applications in nonparametric statistics and image analysis, Khmaladze and Toronjadze [11] ask whether $A_n$ is a consistent estimator for $A$. More precisely, with $|\cdot|$ denoting Lebesgue measure and $\triangle$ denoting symmetric difference of sets, they ask whether we have almost sure convergence

$$|A_n| \to A \qquad \text{as } n \to \infty, \tag{5.1}$$

$$|A \triangle A_n| \to 0 \qquad \text{as } n \to \infty. \tag{5.2}$$

They answer these questions affirmatively only for the case $d = 1$ and comment that for general $d$, (5.2) is not hard to prove when $A$ has Lebesgue null boundary. Using our general results, we can answer these questions affirmatively without any assumptions on the boundary of $A$.

**Theorem 5.1.** *Suppose* $\inf\{\kappa(x) : x \in \Omega\} > 0$. *(5.1) and (5.2) then hold almost surely.*

Note that (5.2) implies (5.1). We keep these results separate for presentational purposes. Actually, the question posed in [11] refers to the almost sure limits analogous to (5.1) and (5.2) for $A_{N_n}$, where $N_n$ is Poisson with parameter $n$ independent of $(X_1, X_2, \ldots)$, but this clearly follows from our result since $N_n \to \infty$ almost surely.

We work toward proving Theorem 5.1. Assume henceforth in this section that $\kappa$ is bounded away from zero on $\Omega$ and set $\Omega_\lambda = \Omega$ for all $\lambda$. For finite $\mathcal{X} \subset \mathbb{R}^d$ and $\mathbf{x} \in \mathcal{X}$, let $\xi(\mathbf{x}; \mathcal{X}, \cdot)$ be the restriction of Lebesgue measure to $\tilde{V}(\mathbf{x}; \mathcal{X})$. Thus, $\xi$ is translation invariant and points do not carry marks; also, $\xi$ has the homogeneity property of order $d$, which says that $\xi(a\mathbf{x}; a\mathcal{X}, aA) = a^d \xi(\mathbf{x}; \mathcal{X}, A)$ for any $a > 0$. Combining this with the consequence (2.2) of translation invariance, we have, for all $\mathbf{x}, \mathcal{X}, A, \lambda$ with $\mathbf{x} \in \Omega$, that

$$\xi_\lambda(\mathbf{x}; \mathcal{X}, A) = \lambda \xi(\mathbf{x}; \mathcal{X}, A) = \lambda |\tilde{V}(\mathbf{x}; \mathcal{X}) \cap A|. \tag{5.3}$$

**Lemma 5.1.** *There is a constant $C$ such that, for $t \geq 1$,*

$$\sup_{n \geq 1} P[\xi_n(\mathbf{X}; \mathcal{X}_{n-1}, \Omega) > t] \leq C \exp(-t/C). \tag{5.4}$$

**Proof.** Let $\mathcal{C}_i, 1 \leq i \leq I$, be a finite collection of infinite open cones in $\mathbb{R}^d$ with angular radius $\pi/12$ and apex at $\mathbf{0}$, with union $\mathbb{R}^d$. For $\mathbf{x} \in \Omega$ and $1 \leq i \leq I$, let $\mathcal{C}_i(\mathbf{x})$ be the translate of $\mathcal{C}_i$ with apex at $\mathbf{x}$. Let $\mathcal{C}_i^+(\mathbf{x})$ be the open cone concentric to $\mathcal{C}_i(\mathbf{x})$ with apex $\mathbf{x}$ and angular radius $\pi/6$. Let $R_{i,n}(\mathbf{x})$ denote the distance from $\mathbf{x}$ to the nearest point in $\mathcal{X}_n \cap \mathcal{C}_i^+(\mathbf{x}) \cap B_{\text{diam}(\mathcal{C}_i(\mathbf{x}) \cap \Omega)}(\mathbf{x})$; if no such point exists, set $R_{i,n}(\mathbf{x}) := \text{diam}(\mathcal{C}_i(\mathbf{x}) \cap \Omega)$. In other words, set

$$R_{i,n}(\mathbf{x}) := \min(\min\{|\mathbf{Y} - \mathbf{x}| : \mathbf{Y} \in \mathcal{X}_n \cap \mathcal{C}_i^+(\mathbf{x})\}, \text{diam}(\mathcal{C}_i(\mathbf{x}) \cap \Omega)),$$

with the convention that $\min(\varnothing) := +\infty$. By elementary geometry, if $\mathbf{Y} \in \mathcal{X}_n \cap \mathcal{C}_i^+(\mathbf{x})$, then $\tilde{V}(\mathbf{x}; \mathcal{X}_n) \cap \mathcal{C}_i^+(\mathbf{x}) \subseteq B_{|\mathbf{Y}-\mathbf{x}|}(\mathbf{x})$. Hence, $V(\mathbf{x}; \mathcal{X}_n) \cap \mathcal{C}_i(\mathbf{x}) \subseteq B_{R_{i,n}}(\mathbf{x})$. Therefore, recalling that $\omega_d$ denotes the volume of the unit ball in $\mathbb{R}^d$, we have, by (5.3), for any $\mathbf{x} \in \Omega$, that

$$\xi_\lambda(\mathbf{x}; \mathcal{X}_{n-1}, \Omega) = \lambda |V(\mathbf{x}; \mathcal{X}_{n-1})| \leq \omega_d \lambda \max_{1 \leq i \leq I} R_{i,n-1}(\mathbf{x})^d. \tag{5.5}$$

Let $\eta := (1/2) \sin(\pi/12)$. Then, $P[R_{i,n}(\mathbf{x}) \geq s] = 0$, unless there exists $\mathbf{y} \in \mathcal{C}_i(\mathbf{x}) \cap \Omega$ with $|\mathbf{y} - \mathbf{x}| = s$. But, in this case, $B_{\eta s}(\frac{\mathbf{x}+\mathbf{y}}{2}) \subseteq \mathcal{C}_i^+(\mathbf{x})$ so that $R_{i,n}(\mathbf{x}) \leq s$, unless $B_{\eta s}(\frac{\mathbf{x}+\mathbf{y}}{2})$ contains no point of $\mathcal{X}_n$. Moreover, since $\frac{\mathbf{x}+\mathbf{y}}{2} \in \Omega$ and $\kappa$ is bounded away from zero on $\Omega$, there is a constant $\delta$, independent of $\mathbf{x}$, such that $\kappa(B_{\eta s}(\frac{\mathbf{x}+\mathbf{y}}{2})) \geq \delta s^d$. Hence, for $1 \leq i \leq I$ and all $u > 0$,

$$P[R_{i,n-1}(\mathbf{x})^d \geq u] \leq (1 - \delta u)^{n-1} \leq \exp(-\delta(n-1)u). \tag{5.6}$$

By (5.5) and (5.6),

$$P[\xi_n(\mathbf{X}; \mathcal{X}_{n-1}, \Omega) > t] \leq \sum_{i=1}^I P[R_{i,n-1}(\mathbf{X})^d > t/(\omega_d n)]$$

$$\leq I \exp(-\delta(n-1)t/(\omega_d n))$$

and this gives us the result. $\square$

**Proof of Theorem 5.1.** Let $\xi_\lambda^*(\mathbf{x}; \mathcal{X}, \cdot)$ be the point mass at $\mathbf{x}$ of size $\xi_\lambda(\mathbf{x}; \mathcal{X}, \Omega)$, as defined at (2.4). Then by (5.3), for Borel $A \subseteq \mathbb{R}^d$, we have

$$\lambda^{-1} \nu_{\lambda,n}^*(A) = \lambda^{-1} \sum_{\mathbf{x} \in A \cap \mathcal{X}_n} \xi_\lambda(\mathbf{x}; \mathcal{X}_n, \Omega) = \sum_{\mathbf{x} \in A \cap \mathcal{X}_n} |V(\mathbf{x}; \mathcal{X}_n)| = |A_n|.$$

Let $\mathcal{C}_i, 1 \leq i \leq I$, be as in the proof of Lemma 5.1. Then, $\mathcal{H}_\lambda \cap \mathcal{C}_i \neq \varnothing$ almost surely for each $i \leq I$; set $R_i(\lambda) := \inf\{|\mathbf{x}| : \mathbf{x} \in \mathcal{H}_\lambda \cap \mathcal{C}_i\}$ and $R(\lambda) := 2 \max_{1 \leq i \leq I} R_i(\lambda)$. The Voronoi cell around $\mathbf{0}$ is unaffected by changes to $\mathcal{H}_\lambda$ outside $B_{R(\lambda)}$ and hence, with this choice of $\xi$, we have $R(\mathbf{0}; \mathcal{H}_\lambda) < \infty$ almost surely for all $\lambda > 0$. Also, $\xi_\infty^{\mathbf{x}}(\mathcal{H}_\lambda, \mathbb{R}^d)$ is the Lebesgue measure of the cell centred at $\mathbf{0}$ in the Voronoi tessellation of $\mathcal{H}_\lambda \cup \{\mathbf{0}\}$. Since $\mathbf{y}$ lies in this cell if and only if $B_{|\mathbf{y}|}(\mathbf{y})$ contains no point of $\mathcal{H}_\lambda$, by Fubini's theorem, we have

$$\mathbb{E}[\xi_\infty^{\mathbf{x}}(\mathcal{H}_\lambda, \mathbb{R}^d)] = \int_{\mathbb{R}^d} P[\mathcal{H}_\lambda \cap B_{|\mathbf{y}|}(\mathbf{y}) = \varnothing] \, d\mathbf{y}$$

$$= \int_{\mathbb{R}^d} \exp(-\lambda \omega_d |\mathbf{y}|^d) \, d\mathbf{y} = 1/\lambda.$$

Set $\lambda(n) = n$ for all $n$. By Lemma 5.1, the measure $\xi$ satisfies the moments condition (2.5). Also, $\xi$ satisfies assumption A2. Hence, setting $f$ to be the indicator function $\mathbf{1}_A$, we can apply Theorem 2.1 to deduce that

$$n^{-1}\nu_{n,n}^*(A) \xrightarrow{L^2} \int_A (1/\kappa(\mathbf{x}))\kappa(\mathbf{x})\,\mathrm{d}\mathbf{x} = |A|,$$

that is, we have (5.1) with $L^2$ convergence.

With $f = \mathbf{1}_A$, using (5.3), we have

$$\lambda^{-1}\langle f, \xi_\lambda(\mathbf{x};\mathcal{X}) - \xi_\lambda^*(\mathbf{x};\mathcal{X})\rangle = \lambda^{-1}(\xi_\lambda(\mathbf{x};\mathcal{X},A) - f(\mathbf{x})\xi_\lambda(\mathbf{x};\mathcal{X},\Omega))$$
$$= \xi(\mathbf{x};\mathcal{X},A) - f(\mathbf{x})\xi(\mathbf{x};\mathcal{X},\Omega)$$
$$= (1 - f(\mathbf{x}))\xi(\mathbf{x};\mathcal{X},A) - f(\mathbf{x})\xi(\mathbf{x};\mathcal{X},\Omega \setminus A)$$

and hence

$$\sum_{\mathbf{x}\in\mathcal{X}} \lambda^{-1}|\langle f,\xi_\lambda(\mathbf{x};\mathcal{X}) - \xi_\lambda^*(\mathbf{x};\mathcal{X})\rangle| = \sum_{\mathbf{x}\in\mathcal{X}\setminus A} \xi(\mathbf{x};\mathcal{X},A) + \sum_{\mathbf{x}\in\mathcal{X}\cap A} \xi(\mathbf{x};\mathcal{X},\Omega\setminus A)$$
$$= \sum_{\mathbf{x}\in\mathcal{X}\setminus A} |V(\mathbf{x};\mathcal{X}) \cap A| + \sum_{\mathbf{x}\in\mathcal{X}\cap A} |V(\mathbf{x};\mathcal{X})\setminus A|$$
$$= \left| A \triangle \bigcup_{\mathbf{x}\in\mathcal{X}\cap A} V(\mathbf{x};\mathcal{X}) \right| = |A \triangle A_n|.$$

Therefore, by applying the conclusion (2.7) of the general result in this particular case with $\lambda(n) = n$, we obtain (5.2) with $L^1$ convergence.

For the almost sure convergence, we demonstrate the condition (2.11) for the present choice of $\xi$. Observe that for $1 \leq i \leq n-1$, the signed measure $\xi(\mathbf{X}_i;\mathcal{X}_{n-1}) - \xi(\mathbf{X}_i;\mathcal{X}_n)$ is, in fact, a non-negative measure, namely the Lebesgue measure on $\tilde{V}(\mathbf{X}_i;\mathcal{X}_{n-1}) \cap \tilde{V}(\mathbf{X}_n;\mathcal{X}_n)$, since this is the region (if any) removed from the Voronoi cell around $\mathbf{X}_i$ due to the addition of an extra point at $\mathbf{X}_n$. Thus, by (5.3), $\xi_\lambda(\mathbf{X}_i;\mathcal{X}_{n-1}) - \xi_\lambda(\mathbf{X}_i;\mathcal{X}_n)$ is $\lambda$ times the same measure. Hence $\xi_\lambda(\mathbf{X}_i;\mathcal{X}_n) - \xi_\lambda(\mathbf{X}_i;\mathcal{X}_{n-1})$ has no positive part and its total variation on $\Omega$ is

$$|\xi_\lambda(\mathbf{X}_i;\mathcal{X}_n) - \xi_\lambda(\mathbf{X}_i;\mathcal{X}_{n-1})|(\Omega) = \lambda|V(\mathbf{X}_i;\mathcal{X}_{n-1}) \cap V(\mathbf{X}_n;\mathcal{X}_n)|.$$

Hence, by (2.8) and (5.3),

$$\tilde{\Delta}_{n,n} = \xi_n(\mathbf{X}_n;\mathcal{X}_n,\Omega) + \sum_{j=1}^{n-1} |\xi_n(\mathbf{X}_j;\mathcal{X}_n) - \xi_n(\mathbf{X}_j;\mathcal{X}_{n-1})|(\Omega)$$
$$= 2n|V(\mathbf{X}_n;\mathcal{X}_n)|$$

and the third moments of this are bounded uniformly in $n$ by Lemma 5.1. Thus, (2.11) holds here with $\beta = 1$ and $p' = 3$, so both (2.6) and (2.7) hold with almost sure convergence. □

### 5.2. Nonparametric regression: the Gamma test

Suppose that $X_i, i \geq 1$, are independent random $d$-vectors with common density $\kappa$. In a nonparametric regression model, consider real-valued random variables $Y_i, 1 \leq i \leq n$, related to random $d$-vectors $X_i, 1 \leq i \leq n$, by the relation

$$Y_i = h(X_i) + e_i, \qquad 1 \leq i \leq n, \tag{5.7}$$

where $h \in C^2(\mathbb{R}^d, \mathbb{R})$ and $(e_i, i \geq 1)$ are independent and identically distributed with mean zero and common variance $\sigma^2$, independent of $(X_i)_{i \geq 1}$. Both the function $h$ and the variance $\sigma^2$ are unknown and $\kappa$ may also be unknown. Often, it is of primary interest to estimate $h$, but here, we are concerned with estimating $\sigma^2$.

Given $k \in \mathbb{N}$ and $n \geq k$, $1 \leq i \leq n$, let $j(i,n,k)$ be the index of the $k$th nearest neighbour of $X_i$ in the sample $\{X_1, \ldots, X_n\}$, that is, the value of $j$ such that $|X_\ell - X_i| \leq |X_j - X_i|$ for precisely $k-1$ values of $\ell \in \{1, 2, \ldots, n\} \setminus \{i\}$. The so-called *Gamma statistic* discussed by Evans and Jones [5] (see also [10]) with parameter $k$ is an estimator $\gamma_{n,k}$ for $\sigma^2$ given by

$$\gamma_{n,k} := \frac{1}{2n} \sum_{i=1}^{n} (Y_{j(i,n,k)} - Y_i)^2. \tag{5.8}$$

For large $n$, we expect $|X_{j(i,n,k)} - X_i|$ to be small, so we approximate to $h(Y_{j(i,n,k)}) - h(Y_i)$ by the first-order Taylor approximation $\nabla h(X_i) \cdot (X_{j(i,n,k)} - X_i)$. Under the proposed model, this approximation gives us

$$\gamma_{n,k} - \sigma^2 \approx \frac{1}{2n} \sum_{i=1}^{n} [(\nabla h(X_i) \cdot (X_{j(i,n,k)} - X_i) + e_{j(i,n,k)} - e_i)^2 - 2\sigma^2]. \tag{5.9}$$

The mean of the last expression is $\frac{1}{2}\mathbb{E}[(\nabla h(X_1) \cdot (X_{j(1,n,k)} - X_1))^2]$. Following equation (2.9) of [5], define

$$A_{n,k} := \frac{\mathbb{E}[(\nabla h(X_1) \cdot (X_{j(1,n,k)} - X_1))^2]}{2\mathbb{E}[|X_{j(1,n,k)} - X_1|^2]}. \tag{5.10}$$

Evans and Jones [5] set $\delta_{k,n} := n^{-1} \sum_{i=1}^{n} |X_{j(i,n,k)} - X_i|^2$ and propose to estimate $\sigma^2$ by the intercept on the $y$-axis of a regression of $y = \gamma_{n,k}$ against $x = \delta_{n,k}$, plotted for $1 \leq k \leq k_0$, with, for example, $k_0 = 20$. They argue heuristically (see the discussion leading up to Theorem 2.1 of [5]) that for large $n$, the value of $A_{n,k}$ should be approximately independent of $k$ and give the slope of the regression line. The following result proves the first of these assertions as an asymptotic statement since the right-hand side of (5.11) below does not depend on $k$. In proving this, we shall give the asymptotic behaviour of the expected value of the expression in (5.9).

We assume throughout this section that $\Omega \subset \mathbb{R}^d$ is bounded and open with $\kappa(\Omega) = 1$, that $|\nabla h|$ is bounded on $\Omega$ and that $\kappa(B_r(x))/r^d$ is bounded away from zero, uniformly over $x \in \Omega, 0 < r \leq 1$. The last condition holds, for example, if $\Omega$ is a finite union of convex sets and the density function $\kappa$ is bounded away from zero on $\Omega$.

**Theorem 5.2.** *As $n \to \infty$,*

$$A_{n,k} \to \frac{\int_\Omega \kappa(x)^{(d-2)/d} |\nabla h(x)|^2 \, dx}{2d \int_\Omega \kappa(x)^{(d-2)/d} \, dx}. \tag{5.11}$$

As $n$ gets larger, one expects the Gamma estimator of $\sigma^2$ to become more sensitive, thus one expects to be able to estimate smaller values of $\sigma^2$. In the next result, we allow the common variance of the $Y_i - h(X_i)$ to get smaller as $n$ increases. More precisely, we modify (5.7) to

$$Y_{i,n} = h(X_i) + n^{-1/d} e_i, \qquad 1 \leq i \leq n. \tag{5.12}$$

We consider an estimator for $\sigma^2$ in this model using a linear regression of just two points arising from $k = 1$ and $k = 2$ (we could similarly consider any other two choices of $k$). Let the random variable $\rho_k$ denote the distance from the origin $\mathbf{0}$ to its $k$th nearest neighbour in the point set $\mathcal{H}_1$. It is known (see [22], equation (19)) that

$$\mathbb{E}[\rho_k^2] = \omega_d^{-2/d} \Gamma(k + (2/d))/\Gamma(k).$$

**Theorem 5.3.** *Let $\gamma_{n,k}$ be given by (5.8), with $Y_i = Y_{i,n}$ given by (5.12). Then,*

$$n^{2/d} \left( \frac{\gamma_{n,2} \mathbb{E}\rho_1^2 - \gamma_{n,1} \mathbb{E}\rho_2^2}{\mathbb{E}\rho_1^2 - \mathbb{E}\rho_2^2} \right) \xrightarrow{L^1} \sigma^2. \tag{5.13}$$

This result shows that the left-hand side of (5.13) (which is the intercept in linear regression of $n^{2/d} \gamma_{n,k}$ against $\mathbb{E}\rho_k$, based on just two values of $k$) is a consistent estimator of $\sigma^2$.

In proving Theorems 5.2 and 5.3, we use the following notation. For locally finite $\mathcal{X} \subset \mathbb{R}^d$ and for $x \in \mathbb{R}^d$, let $N_\mathcal{X}^k(x)$ be the $k$th nearest neighbour of $x$ in the set $\mathcal{X} \setminus \{x\}$, making the choice according to an arbitrary rule in the event of ties and taking $N_\mathcal{X}^k(x) = x$ if $\mathcal{X} \setminus \{x\}$ has fewer than $k$ elements.

**Lemma 5.2.** *Suppose $\mathbf{b} \in \mathbb{R}^d$ and $\lambda > 0$. Then,*

$$\mathbb{E}[|\mathbf{b} \cdot N_{\mathcal{H}_\lambda}^k(\mathbf{0})|^2] = d^{-1} \lambda^{-2/d} |\mathbf{b}|^2 \mathbb{E}[\rho_k^2]. \tag{5.14}$$

**Proof.** Let $\Theta = (\theta_1, \ldots, \theta_d)$ be uniformly distributed over the unit sphere in $\mathbb{R}^d$, independent of $\rho_k$. Since $\sum_{i=1}^{d} \theta_i^2 = |\Theta|^2 = 1$, taking expectations, we have $\mathbb{E}[\theta_i^2] = 1/d$ for each $i$. Given $\mathbf{b}$, we have, for some unit vector $\mathbf{e}$, that $\mathbf{b} = |\mathbf{b}|\mathbf{e}$ and hence

$$\mathbb{E}[(\mathbf{b} \cdot \Theta)^2] = |\mathbf{b}|^2 \mathbb{E}[|\mathbf{e} \cdot \Theta|^2] = d^{-1}|\mathbf{b}|^2, \tag{5.15}$$

where the last equality follows because, by rotational symmetry, the distribution of $|\mathbf{e} \cdot \Theta|^2$ is the same for all unit vectors $\mathbf{e}$ and its expectation is $1/d$ whenever $\mathbf{e}$ is one of the unit coordinate vectors.

By the distributional rotational symmetry of the homogeneous Poisson process and the fact that $\mathcal{H}_\lambda$ has the same distribution as $\lambda^{-1/d} \mathcal{H}_1$ for any $\lambda > 0$ (by the Mapping Theorem [12]), we have $N_{\mathcal{H}_\lambda}^k(\mathbf{0}) \stackrel{\mathcal{D}}{=} \lambda^{-1/d} \rho_k \Theta$. Taking expectations and using (5.15), we have (5.14). □

**Proof of Theorem 5.2.** Let $\xi(x; \mathcal{X})$ be a point mass at $x$ of size $(\nabla h(x) \cdot (N_{\mathcal{X}}^k - x))^2$ and let $\Omega_\lambda = \Omega$ for all $\lambda$. Then, for $x \in \Omega$, $\xi_\lambda(x; \mathcal{X})$ is a point mass at $x$ of size $\lambda^{2/d}(\nabla h(x) \cdot (N_{\mathcal{X}}^k(x) - x))^2$. By assumption, $|\nabla h|$ is bounded on $\Omega$, so there is a constant $C$ such that, for $x \in \Omega$,

$$|\nabla h(x) \cdot (N_{\mathcal{X}}^k(x) - x)|^2 \leq C|N_{\mathcal{X}}^k(x) - x|^2$$

and hence, for $x \in \Omega$,

$$P[\xi_n(x; \mathcal{X}_{n-1}, \Omega) > t] \leq P[|N_{\mathcal{X}_{n-1}}^k(x) - x|^2 > t/(Cn^{2/d})].$$

By assumption, $\kappa(B_r(x))/r^d$ is bounded away from zero on $0 < r < \mathrm{diam}(\Omega)$, so there are constants $C', C''$ such that for $x \in \Omega$, $n \geq 2k$ and $(t/(Cn^{2/d}))^{1/2} \leq \mathrm{diam}(\Omega)$,

$$P[\xi_n(x; \mathcal{X}_{n-1}, \Omega) > t] \leq \sum_{j=0}^{k-1} \binom{n-1}{j} \left(\frac{t^{d/2}}{C'n}\right)^j \left(1 - \left(\frac{t^{d/2}}{C'n}\right)\right)^{n-1-j}$$

$$\leq \exp(-t^{d/2}/C'')$$

and this bound also holds for $(t/(Cn^{2/d}))^{1/2} > \mathrm{diam}(\Omega)$ since, in that case, the probability is zero. It follows that (2.5) holds here for any $p$. Therefore, we may apply Theorem 2.1, here taking $f \equiv 1$ and $\Omega_\lambda = \Omega$ for all $\lambda$, followed by Lemma 5.2, to obtain

$$n^{2/d} \mathbb{E}[(\nabla h(X_1) \cdot (X_{j(i,n,k)}) - X_1)^2] = \mathbb{E}[\xi_n(X_1; \mathcal{X}_n, \Omega)] = n^{-1} \mathbb{E}[\langle f, \nu_{n,n} \rangle]$$

$$\to \int_\Omega \mathbb{E}[(\nabla h(x) \cdot N_{\mathcal{H}_{\kappa(x)}}^k(\mathbf{0}))^2] \kappa(x) \, dx$$

$$= d^{-1} \mathbb{E}[\rho_k^2] \int_\Omega \kappa(x)^{1-(2/d)} |\nabla h(x)|^2 \, dx. \tag{5.16}$$

As for the denominator of the expression (5.10) for $A_{n,k}$, we have (see [6, 18, 22]) that

$$n^{2/d}\mathbb{E}[|X_{j(i,n,k)} - X_1|^2] = n^{-1}\mathbb{E}\sum_{j=1}^{n}(n^{1/d}|X_{j(i,n,k)} - X_i|)^2$$

$$\to \int_{\Omega}(\kappa(x)^{-2/d}\mathbb{E}[\rho_k^2])\kappa(x)\,\mathrm{d}x$$

and combining this with the limiting expression (5.16) for the numerator, we obtain (5.11). □

**Proof of Theorem 5.3.** We give just a sketch. We must now consider marked points. The mark space is $\mathcal{M} = \mathbb{R}$ and the real-valued mark $T(x)$ attached to a point $x \in \mathbb{R}^d$ is assumed to have the common distribution of $e_1, e_2, \ldots$ (with mean zero and variance $\sigma^2$). Given $k$ and given $\mathbf{x} = (x, T(x))$ with $\mathbf{x} \in \mathcal{X}$ and $\mathcal{X}$ a finite subset of $\mathbb{R}^d \times \mathcal{M}$, we set

$$\xi(\mathbf{x}; \mathcal{X}) = (\nabla h(x) \cdot (N_{\mathcal{X}}^k(x) - x) + T(N_{\mathcal{X}}^k(x)) - T(x))^2.$$

This can be shown to satisfy the conditions of Theorem 2.1. We set $e_i = T(X_i)$ and $\mathbf{X}_i = (X_i, e_i)$ and $\mathcal{X}_n = \{\mathbf{X}_1, \ldots, \mathbf{X}_n\}$. By first-order Taylor approximation,

$$n^{2/d}\gamma_{n,k} \approx \frac{n^{(2/d)-1}}{2}\sum_{i=1}^{n}(\nabla h(X_i) \cdot (X_{j(i,n,k)} - X_i) + n^{-1/d}(e_{j(i,n,k)} - e_i))^2$$

$$= \frac{1}{2n}\sum_{i=1}^{n}(n^{1/d}(\nabla h(X_i) \cdot (X_{j(i,n,k)} - X_i)) + (e_{j(i,n,k)} - e_i))^2$$

$$= \frac{1}{2n}\sum_{i=1}^{n}\xi_n(X_i, \mathcal{X}_n) = \frac{1}{2n}\langle f, \nu_{n,n}\rangle,$$

where we put $f \equiv 1$. By Theorem 2.1, this converges in $L^1$ to the limit

$$\frac{1}{2}\int_{\Omega}\mathbb{E}[(\nabla h(x) \cdot N_{\mathcal{H}_{\kappa(x)}} + T(N_{\mathcal{H}_{\kappa(x)}}) - T(\mathbf{0}))^2]\kappa(x)\,\mathrm{d}x$$

$$= \sigma^2 + \frac{\mathbb{E}[\rho_k^2]}{2d}\int_{\Omega}|\nabla h(x)|^2\kappa(x)^{1-(2/d)}\,\mathrm{d}x,$$

where we have used Lemma 5.2 and the fact that the marks $T$ have mean zero, variance $\sigma^2$ and are independent of each other and of the point process $\mathcal{H}_{\kappa(x)}$.

It then follows that we have the $L^1$ convergence

$$n^{2/d}(\gamma_{n,2}\mathbb{E}\rho_1^2 - \gamma_{n,1}\mathbb{E}\rho_2^2) \to (\mathbb{E}\rho_1^2 - \mathbb{E}\rho_2^2)\sigma^2$$

and this implies (5.13). □